\documentstyle{amsppt}
%%\redefine\publogo{}

\redefine\widetilde{\mathaccent"0365}
\redefine\leq{\leqslant}
\redefine\geq{\geqslant}
\def\myone{\hbox{\bf 1}}

\def\arctg{\operatorname{arctg}}

\def\cond{\operatorname{cond}}

\def\tg{\operatorname{tg}}
\def\th{\operatorname{th}}
\define\DD{\text{D}}
\def\Arcth{\operatorname{Arcth}}
%%proizvodnaja:
\define\pd#1#2{\dfrac{\partial #1}{\partial#2}}
\define\sumpr{\sideset \and^{\prime} \to\sum}

\NoRunningHeads
\topmatter

\title Approximate construction of rational approximations and the
effect of error autocorrection. Applications.\endtitle

\author G. L. Litvinov\endauthor
\thanks{Published in: Russian Journal of Mathematical
Physics, vol.1, No. 3, 1994.}\endthanks
%%\endtopmatter

\abstract Several construction methods for rational approximations to functions of
one real variable are described in the present paper; the computational
results
that characterize the comparative accuracy of these methods are presented; an
effect of error autocorrection is considered.
This effect occurs in efficient methods of rational approximation (e.g., Pad\' e
approximations, linear and nonlinear Pad\' e--Chebyshev approximations)
where very significant errors in the coefficients do not affect the
accuracy of the approximation. The matter of import is that the errors in the
numerator and the denominator of a fractional rational approximant compensate
each other. This effect is related to the fact that the errors in the
coefficients of a rational approximant are not distributed in an arbitrary
way but form the coefficients of a new approximant to the approximated
function. Understanding of the error autocorrection mechanism allows to
decrease this error by varying the approximation procedure depending on the
form of the approximant. Some applications are described in the paper. In
particular, a method of implementation of basic calculations on decimal
computers that uses the
technique of rational approximations is described in the Appendix.

To a considerable extent the paper is a survey and the
exposition is as elementary as possible.\endabstract

\toc

\specialhead  {}\S1. Introduction\page{2}\endspecialhead
\specialhead {}\S2. Best approximants\page{4}\endspecialhead
\specialhead {}\S3. Construction methods for best approximants
\page{6}\endspecialhead
\specialhead {}\S4. The role of approximate methods
and an estimate of the quality of approximation\page{8}\endspecialhead
\specialhead {}\S5. Chebyshev polynomials and polynomial approximations
\page{10}\endspecialhead
\specialhead {}\S6. Ill-conditioned problems and rational approximations
\page{13}\endspecialhead
\specialhead {}\S7. The effect of error autocorrection
\page{15}\endspecialhead
\specialhead {}\S8. Pad\'e approximations\page{17}\endspecialhead
\specialhead {}\S9. Linear Pad\'e--Chebyshev approximations and the
PADE program\page{19}\endspecialhead
\specialhead {}\S10. The PADE program. Analysis of the algorithm
\page{22}\endspecialhead
\specialhead {}\S11. The ``cross--multiplied''
linear Pad\'e--Chebyshev approximation scheme\page{26}\endspecialhead
\specialhead {}\S12. Nonlinear Pad\'e--Chebyshev approximations
\page{27}\endspecialhead
\specialhead {}\S13. Applications of the computer algebra system
REDUCE to the construction of rational approximants\page{28}\endspecialhead
\specialhead {}\S14. The effect of error autocorrection for nonlinear
Pad\'e--Cheby\-shev approximations\page{29}\endspecialhead
\specialhead {}\S15. Small deformations of approximated functions and
acceleration of convergence of series\page{32}\endspecialhead
\specialhead {}\S16. Applications to computer calculation\page{33}\endspecialhead
\specialhead {}\S17. Nonlinear models and rational approximants
\page{34}\endspecialhead
\specialhead {} APPENDIX. A method of implementation of basic
calculations on decimal computers\page{35}\endspecialhead
 \subhead 1. Introduction\page{35}\endsubhead
 \subhead 2. Floating point arithmetic system\page{36}\endsubhead
 \subhead 3. The design of computation for elementary functions
\page{36}\endsubhead
 \subhead 4. Algorithms\page{37}\endsubhead
 \subsubhead 4.1. Calculation of logarithms\page{37}\endsubsubhead
 \subsubhead 4.2. Calculation of exponentials\page{38}\endsubsubhead
 \subsubhead 4.3. Calculation of $\sin x$ and $\cos x$\page{38}\endsubsubhead
 \subsubhead 4.4. Calculation of $tg x$\page{39}\endsubsubhead
 \subsubhead 4.5. Calculation of $arctg x$\page{39}\endsubsubhead
 \subsubhead 4.6. Calculation of $\arcsin x$\page{39}\endsubsubhead
 \subhead 5. Coefficients\page{39}\endsubhead
 \subhead 6. Analysis of the algorithms\page{39}\endsubhead
 \subhead 7. Implementation of algorithms for calculating elementary
functions\page{40}\endsubhead
\specialhead {} References\page{43}\endspecialhead
\endtoc
\endtopmatter

\document

\leftskip = 2 in
{\eightpoint
\it Whenever he has some money to spare, he goes to a shop and buys some
kind of useful book. Once he bought a book that was entitled ``Inverse
trigo\-nomet\-ri\-cal functions and Chebyshev polynomials''.

{\rom{N.N.Nosov}} ``{\it Happy family}''. Moscow, 1975, p.91}

\leftskip = 0 in

\head \S1. Introduction\endhead

The author came across the phenomenon of error autocorrection at the end of
seventies while developing nonstandard algorithms for computing elementary
functions on small computers. It was required to construct rational
approximants of the form
$$
R(x)={a_0+a_1x+a_2x^2+\dots +a_nx^n\over{b_0+b_1x+b_2x^2+\dots +b_mx^m}}
\tag1
$$
to certain functions of one variable $x$ defined on finite segments of the
real line. For this purpose a simple method (described in~\cite{1} and
below) was used: the method allows to determine the family of
coefficients $a_i$, $b_j$ of the approximant~(1) as the solution of a
certain system of
linear algebraic equations. These systems
turned out to be ill conditioned, i.e., the problem of
determining the coefficients of the approximant is, generally speaking,
ill-posed in the sense of~\cite{2}. Nevertheless, the method ensures a
paradoxically high quality of the obtained approximants whose errors are
close to the best possible~\cite{1}.

For example, for the function $\cos x$ the approximant of the form~(1) on
the segment $[-\pi /4,\pi /4]$ obtained by the method mentioned above for
$m=4$, $n=6$ has the relative error equal to
$0.55\cdot 10^{-13}$, and the best
possible relative error is $0.46\cdot 10^{-13}$ \cite{3}. The corresponding
system of linear algebraic equations has the condition number of order
$10^9$. Thus we risk losing 9 accurate decimal digits in the
solution of calculation errors. Computer experiments
show that this is a serious risk. The method mentioned above
was implemented as a Fortran program. The calculations were
carried out with double precision (16 decimal positions) by means of
ICL--4--50 and ES--1045 computers. These computers are very similar in their
architecture, but when passing from one computer to another the system of
linear equations and the computational process are perturbed because of
calculation errors, including round-off errors. As a result, the
coefficients of the approximant mentioned above to the function $\cos x$
experience a perturbation already at the
sixth--ninth decimal digits. But the error in the rational approximant
itself remains
invariant and is $0.4\cdot 10^{-13}$ for the absolute error and $0.55\cdot
10^{-13}$ for the relative error. The same thing happens for approximants
of the
form~(1) to the function $\arctg x$ on
the segment [-1,1] obtained by the method mentioned above
for $m=8$, $n=9$ the relative error is $0.5\cdot 10^{-11}$
and does not change while passing from ICL--4--50 to ES--1045 although the
corresponding system of linear equations has the condition number of
order $10^{11}$, and the coefficients of the approximant experience a
perturbation with relative error of order~$10^{-4}$.

Thus the errors in the numerator and the denominator of a rational approximant
compensate each other. The effect of error autocorrection is connected with
the fact that the errors in the coefficients of a rational approximant are not
distributed in an arbitrary way, but form the coefficients of a new
approximant to the approximated function. It can be easily understood that the
standard methods of interval arithmetic (see, for example,~\cite{54})
do not allow to take into account this effect
and, as a result, to estimate the error in the rational approximant
accurately.

Note that the application of standard procedures known
in the theory of ill-posed problems results in this case in losses
in accuracy. For example, if one applies the regularization method,
two thirds of the accurate figures are lost~\cite{4};
in addition, the amount of calculations increases rapidly. The
matter of import is that the exact solution of the system of equations
in the
present case is not the ultimate goal; the aim is to construct an
appro\-ximant which is precise enough. This approach allows to
``rehabilitate'' (i.e., to justify) and to simplify a number of
algorithms intended for the construction of the approximant, to obtain (without
additional transforms) approximants in a form which is convenient for
applications.

Professor Yudell L.~Luke kindly drew the author's attention to his
papers~\cite{5, 6} where the effect of error autocorrection for the classical
Pad\' e approximants was revealed and was explained at a heuristic level.
The method mentioned above leads to the linear
Pad\' e--Chebyshev approximants if the calculation errors are
ignored.

In the present paper, using heuristic arguments and the results of computer
experiments, the error autocorrection mechanism is considered for quite a
general situation (linear methods for the construction of
rational approximants,
nonlinear generalized Pad\' e approximations). The efficiency of the
construction algorithms used for rational approximants seems to be
due to the error autocorrection effect (at least in the case when the number
of coefficients is large enough).

Our new understanding
of the error autocorrection mechanism allows us, to some
extent, to control calculation errors by changing
the construction procedure depending on the form of the approximant.

In the paper the construction algorithm for linear Pad\' e--Chebyshev
approximants is considered and the corresponding program is
briefly described (see~\cite{7}).
It is shown that the appearance of a control
parameter allowing to take into account the error autocorrection mechanism
ensures the decrease of the calculation errors in some cases.
Results of computer calculations that characterize the
possibilities of the program and the quality of the approximants obtained
as compared to the best ones are presented. Some other (linear and
nonlinear) construction methods for
rational approximants are described. Construction methods for linear and
nonlinear Pad\' e--Chebyshev approximants involving the computer
algebra system REDUCE (see~\cite{8}) are also briefly described.
Computation results characterizing the comparative precision of these
methods are given. With regard to the error autocorrection phenomenon the
effect described in~\cite{9} and connected with the fact that a small
variation of an approximated function can lead to a sharp decrease in
accuracy of
the Pad\' e--Chebyshev approximants is analyzed. Some applications are
indicated. In particular,
a method of implementation of basic calculations on decimal computers that
uses the technique of rational approximations is described in the
Appendix.

To a considerable extent the paper is a survey and the exposition is as
elementary as possible. In the survey part of the paper we tried to
present the required information clearly and consistently, to make
it self-contained. But this part does not claim to
be complete: the number of papers concerning rational
approximations theory and its applications in numerical analysis (including
computer calculation of functions, numerical solving of equations,
acceleration of convergence of series, and quadratures), in theoretical and
experimental physics (including quantum field theory, scattering theory,
nuclear and neutron physics), in the theory and
practice of experimental data processing , in mechanics,
in control theory, and other branches is much too
vast; see, in particular, the reviews and reference handbooks~\cite{3,
10--16}.

The author is grateful to Yudell L.~Luke for stimulating conversations
and valuable
instructions. The author also wishes to express his thanks to I.~A.~Andreeva,
A.~Ya.~Rodionov and V.~N.~Fridman who participated in the programming and
organization of computer experiments. This paper would not have been
written without
their help. A preliminary version of the paper was published in~\cite{56}.

\head\S2. Best approximants\endhead

We shall need some information and results pertaining to ideas of
P.~L. Chebyshev, see~\cite{17}. Let $[A,B]$ be a real line segment
(i.e., $[A,B]$ is the set of all real numbers $x$ such that $
A\leq x\leq B$) and $f(x)$ be a continuous function defined on this
segment. Consider the {\it absolute error function} of the
approximant of the form
~(1) to the function $f(x)$, i.e., the quantity
$$
\Delta (x)=f(x)-R(x),\tag 2
$$
and the {\it absolute error} of this approximant, i.e., the quantity
$$
\Delta =\max_{A\leq x\leq B}\vert\Delta (x)\vert =\max_{A\leq x\leq B}\vert
f(x)-R(x)\vert.\tag 3
$$

A classical problem of approximation theory is to determine, for fixed degrees
$m$ and $n$ in~(1), the coefficients in the numerator and the denominator of
expression~(1) so that~(3) is the smallest possible. The corresponding
approximant is called {\it best} (with respect to the absolute error). An
important role is played by the following result.

\proclaim{Generalized de la Vall\'ee--Poussin theorem~\cite{17}}
If the polynomials
$$
\aligned
\widetilde{P}(x)&=\tilde a_0+\tilde a_1x+\dots +\tilde a_{n-\nu}x^{n-\nu},\\
\widetilde{Q}(x)&=\tilde b_0+\tilde b_1x+\dots +\tilde b_{m-\mu}x^{m-\mu},
\endaligned
\tag 4
$$
where $0\leq\mu\leq m$, $0\leq\nu\leq n$, $b_{m-\mu}\neq 0$, have no
common divisor (i.e., the fraction
$\widetilde{R}(x)=\widetilde{P}(x)/\widetilde{Q}(x)$ is irreducible), the
expression $\widetilde{R}(x)=\widetilde{P}(x)/\widetilde{Q}(x)$ is finite on
the segment $[A,B]$, and at successive points $x_1<x_2<\dots <x_N$
of the segment $[A,B]$ the error function
$\tilde\Delta (x)=f(x)-\widetilde{R}(x)$ of the approximant $\widetilde{R}$
takes nonzero values $\lambda_1, -\lambda_2,\dots ,(-1)^N\lambda_N$
with alternating signs (so that the numbers $\lambda_i$ are either all
positive or all negative), $N=m+n+2-d$, where $d$ is the smallest of the
numbers $\mu$, $\nu$, then the error $\Delta$ of any approximant of the
form~(1) satisfies the inequality
$$
\Delta\geq\lambda =\min\big\{\vert\lambda_1\vert,\vert\lambda_2\vert,\dots,
\vert\lambda_N\vert\big\}.\tag5
$$
\endproclaim

\demo{Proof} Suppose that there exists an approximant $R(x)$ of the form~(1)
for which the inequality~(5) is not satisfied. Consider the difference
$$
\varepsilon (x)=R(x)-\widetilde{R}(x)=\big[f(x)-\widetilde{R}(x)\big]
-\big[f(x)-R(x)\big]=\widetilde{\Delta}(x)-\Delta (x).
$$
From our assumption it follows that the numbers $\varepsilon (x_1),
\varepsilon (x_2),\dots ,\varepsilon (x_N)$ differ from zero and have
alternate signs. And this, in its turn, by virtue of the continuity of the
function $\varepsilon (x)$ on the segment $[A,B]$ implies that
the function $\varepsilon (x)$ has at least $N-1=m+n+1-d$ zeros inside the
segment $[A,B]$. On the other hand, the definition of the function
$\varepsilon (x)$ implies the equality $\varepsilon (x)=U(x)/V(x)$,
where $U(x)$ and
$V(x)$ are polynomials and the degree of $U(x)$ does not
exceed $m+n-d$. So the function
$\varepsilon (x)$ cannot have more than $N-2=m+n-d$ zeros. This
contradiction proves the theorem.
\enddemo

The quantity $d$ which is mentioned in the theorem is called the
{\it defect of
the approximant} $\widetilde{R}(x)$; in practice usually $d=\mu=\nu=0$. The
generalized de la Vall\'ee--Poussin theorem gives us a sufficient
condition for the approximant
$\widetilde{R}(x)=\widetilde{P}(x)/\widetilde{Q}(x)$, where
$\widetilde{P}(x)$
and $\widetilde{Q}(x)$ are polynomials of the form~(4), to be best. The
points $x_1<x_2<\dots <x_N$ of the segment $[A,B]$ are called {\it
Chebyshev alternation points} for the approximant $\widetilde{R}(x)$
if the error function $\tilde\Delta (x)=f(x)-\widetilde{R}(x)$ at
these points has values which coincide with the
absolute error $\tilde\Delta$ of the approximant $\widetilde{R}$ in
absolute value and are alternate in sign.
In other words, at the points $x_1,\dots ,x_N$ the error function
$\tilde\Delta (x)$ has
extrema with alternating signs which coincide with each
other in absolute value. From the generalized de la Vall\'ee--Poussin
theorem, it follows that the presence of Chebyshev alternation points is
sufficient for the approximant $\widetilde{R}(x)$ to be best.

\proclaim{Chebyshev theorem} The presence of Chebyshev alternation points
is a necessary and sufficient condition under which the approximant is best.
Such an approximant exists and is unique if two
fractions that coincide after cancellation are not regarded as different.
\endproclaim

A comparatively simple proof is given in~\cite{17}. Note that P.~L.~Chebyshev
and Vall\'ee--Poussin considered the case of polynomial approximants. The
general case was first considered by N.~I.~Akhiezer, the
results mentioned being valid also in the case when the expression
$\Delta_\rho (x)=(f(x)-R(x))/\rho (x)$,
where the weight $\rho$ is
nonzero, is taken for the error function; if the weight satisfies certain
additional conditions, then the segment $[A,B]$ need not be assumed
finite~\cite{17}. Note, that for $\rho (x)\equiv 1$ we obtain the
absolute error~(2); and if $f(x)$ has no zeros on the
segment $[A,B]$, then for $\rho (x)=f(x)$ we shall obtain the {\it relative
error function}
$$
\delta(x)={\Delta(x)\over f(x)}={f(x)-R(x)\over f(x)}=1-{R(x)\over f(x)}.
\tag6
$$
Correspondingly, the quantity
$$
\delta =\max_{A\leq x\leq B}\vert\delta(x)\vert
\tag7
$$
is the relative error, and one can speak of the {\it best approximants
with respect to the relative error}.

Suppose that the segment $[A,B]$ is symmetric with respect to zero, i.e.,
$A=-B$. If the function $f(x)$ is even, then it is not difficult to verify
that all its best rational approximants on this segment (in the sense
of the absolute error or of the relative one) are also even functions, so that
one can immediately look for them in the form $R(x^2)=P(x^2)/Q(x^2)$, where $
P$ and $Q$ are polynomials. If the function $f(x)$ is odd, then its best
approximants are also odd functions and one can immediately look for
them in the form $xR(x^2)=xP(x^2)/Q(x^2)$, where $P$ and $Q$ are polynomials.
One can speak of the best approximants with respect to the relative error, if
an odd function $f(x)$ is zero only for $x=0$, is
continuously differentiable, and $f'(0)\neq 0$. In this case $f(x)$ can be
represented in the form $x\varphi (x)$, where $\varphi (x)$ is a continuous
even function that never equals zero. Then describing rational approximants
to the function $f(x)$ with best relative error reduces to solving
the same problem for the even function $\varphi (x)$; indeed,
$$
{f(x)-xR(x^2)\over f(x)}={x\varphi(x)-xR(x^2)\over x\varphi(x)}=
{\varphi(x)-R(x^2)\over \varphi(x)}.
$$

\head\S3. Construction methods for best approximants\endhead

Suppose that a rational approximant of the form~(1) is the best
approximant to a continuous function $f(x)$ on the segment $[A,B]$. For
simplicity, further we shall assume that the defect is zero.
Let $x_1,x_2,\dots ,x_{m+n+2}$ be the Chebyshev alternation points. Then the
error function $\Delta_\rho (x)$ corresponding to the weight $\rho$ (see
above) satisfies the following system of equalities:
$$
\Delta_\rho(x_k)=(-1)^k\lambda,
\tag8
$$
where $\vert\lambda\vert=\Delta_\rho =\max_{A\leq x\leq
B}\vert\Delta_\rho (x)\vert$; $k=1,\dots ,m+n+2$.
For fixed values of $x_1,\dots ,x_{m+n+2}$,
relations~(8) can be regarded as a system of $m+n+2$ equations with respect
to the unknowns $a_i$, $b_j$, $\lambda$, where $i=0,\dots ,n$, $j=0,\dots ,m$.
Since one can multiply the numerator and the denominator of the fraction
$R(x)$ by the same number, we see that one more condition, for example, $b_0
=1$, can be added to system~(8), so that the number of equations coincides
with the number of unknowns. The iteration method of computation of
coefficients in the approximant $R(x)$ (suggested by A.~Ya.~Remez
(see~\cite{18}) for polynomial approximants and generalized to the general
case) is based on this idea. Different versions of the generalized Remez
method were considered in many papers; see, for example~\cite{3, 12, 20--27}.

The approximant is constructed as follows. On the
first step, the initial approximations
$x_1<x_2<\dots <x_{m+n+2}$ to the Chebyshev alternation points are chosen
on the segment $[A,B]$
and the system of equations~(8) is solved. As a result we obtain some rational
approximant $R_1(x)$ with error function
$\Delta_1(x)=(f(x)-R_1(x))/\rho (x)$. For this function the extremum
points are found, and the information obtained is used to modify the set
$\{ x_1,\dots ,x_{m+n+1}\}$. Then the procedure is repeated anew, a new
approximant $R_2(x)$ is obtained, and so on.

Taking into account the fact that $\Delta_\rho (x)=(f(x)-R(x))/\rho (x)$
and $R(x)$ has the form~(1), system~(8) can be rewritten in the form
$$
f(x_k)-{\sum_{i=0}^n a_i(x_k)^i\over\sum_{j=0}^m b_j(x_k)^j}=(-1)^k
\rho (x)\cdot\lambda,
$$
whence, as the result of elementary transformations, we get the system of
equations
$$
\sum_{i=0}^n a_i(x_k)^i+\mu_k(\lambda)\sum_{j=0}^m b_j(x_k)^j=0,
\tag9
$$
where $\mu_k(\lambda)=(-1)^k\rho (x_k)\lambda -f(x_k)$,
$k=1,2,\dots ,m+n+2$.

Note that for a fixed value of $\lambda$ (as well as for the alternation
points\newline
 $x_1,\dots ,x_{m+n+2}$) the coefficients $a_i$, $b_j$ of the
approximant satisfy the system of linear homogeneous algebraic equations~(9).
But $\lambda$ must also be determined; this transforms~(9)
into a nonlinear system of equations which is rather difficult to solve. The
case when it is necessary to find the polynomial approximant, i.e., the case
$m=0$, is an exception to what was just noted. In this case the system~(9)
becomes linear.

The solution of nonlinear system of equations~(9) is usually reduced to
the iterated solution of
systems of linear equations. The following method is comparatively
popular (see, for example,~\cite{3, 12, 21, 22, 25, 28}) and was used to
compile the well-known tables of rational approximants to elementary
and special functions~\cite{3}. Let $b_0=1$ (normalization); then~(9)
takes the following form
$$
\sum_{i=0}^n a_i(x_k)^i+\mu_k(\lambda)\sum_{j=1}^m b_j(x_k)^j
+(-1)^k\rho(x_k)\lambda=f(x_k).
\tag"{($9'$)}"
$$
Substituting a fixed number $\lambda_0$ for $\lambda$ in the
nonlinear terms of system~(9'), we get the linear system
$$
\sum_{i=0}^n a_i(x_k)^i+\mu_k(\lambda_0)\sum_{j=1}^m b_j(x_k)^j
+(-1)^k\rho(x_k)\lambda=f(x_k).\tag10
$$

The iteration process is applied to the initial collection of values of the
critical points $\{ x_k\}$, i.e., of initial approximations to the Chebyshev
alternation points, and to the given value $\lambda_0$. First, from~(10)
one determines the new value of $\lambda$ and substitutes it for
$\lambda$ in the nonlinear terms of equation~(10); then the system of
equations~(10) is solved again, and the next value of $\lambda$ is
determined, and so on. As a result a new value of $\lambda$ and the
collection of the coefficients $a_i$, $b_j$ are defined. The next step is
to determine a new collection of critical points $\{x_k\}$ as extremum
points of the error function for the approximant obtained on the previous
step. Both steps form one cycle of an iteration process. The calculation is
finished when the value of $\lambda$ with precision given in advance
coincides in absolute value with the maximal value of the error function.
A complete text of the corresponding Algol program is given in~\cite{22}.

Unfortunately the iteration process described above can be nonconvergent
even in the case when the initial approximation differs from the solution of
the problem infinitesimally; see~\cite{28}. For some versions of the
Remez method it is proved that the iteration process converges if the
initial approximation is sufficiently good, see~\cite{20, 23--25, 29, 12}.
Nevertheless, in each particular case it is
often difficult to indicate {\it a priori} (i.e., before the start of
calculations) the initial approximant that ensures the
convergence of the iteration process, and for a given initial approximation it
is difficult to verify whether the conditions which ensure the
convergence are satisfied. One of the methods which is applied in practice
is to construct, at first, the best polynomial approximant of degree
$m+n$ (in this case no difficulties arise); next, using the Chebyshev
alternation points of this approximant as the initial collection of critical
points for the iteration process one constructs the best approximant having
the form of a polynomial of degree $m+n-1$ divided by a linear function.
Finally, in the same manner, the degree of the numerator is successively
reduced and the degree of the denominator is successively raised till an
approximant of the required form~(1) is obtained, see~\cite{3, 12}.

Together with  iteration methods for constructing the best rational
approximants, methods of linear and convex
programming are used, see~\cite{18, 30}. Iteration methods, as a rule, are
more efficient~\cite{27}, but cannot be generalized directly to
the case of functions of several variables.

\head\S4. The role of approximate methods
and an estimate of the quality of approximation\endhead

The construction algorithms for the best rational approximants are
comparatively complicated, so simpler methods that give an approximate
solution of the problem are used on a large scale, see, for
example,~\cite{1, 5, 7--9, 11--15, 24, 25, 31--37}. Below we
describe methods which are easily implemented, use
comparatively little computation time and yield
approximants that are close to best. Such an approximant can be used as an
initial approximant for an iteration algorithm which gives the exact result.
The approximant that is best in the sense of the absolute error is not
necessary best
in the sense of the relative error. It is usually important in practice for
both the absolute error and the relative one to be small. So rather than the
best approximants, the approximants constructed by means of an
approximate method and having appropriate absolute and relative errors are
often more convenient. Finally, one can also apply methods giving an
approximate solution of the rational approximation problem to those
cases when the information about an approximated function is incomplete (for
example, there are known values of a function only for a finite number
of the argument values, or there are known only the first terms of the
function expansion in a series, or the initial information contains an
error, and so on).

The generalized de la Vall\'ee--Poussin theorem (see~\S2 above) allows to
estimate the proximity of an approximate solution of the
approximation problem to the best approximant even in the case when this
best approximant is unknown.

For example, suppose we want to estimate the proximity of a given
approximant of the form~(1) to an approximant of the same form
with best absolute error to a given
function $f(x)$. Suppose for simplicity that
the defect of the best approximant is zero (in practice this
condition usually holds). Then, by virtue of Chebyshev's
alternation theorem, in the case when the given approximant $R(x)$ is
sufficiently close to the best one, at the
successive points $x_1<\dots <x_{m+n+2}$ belonging to the interval where the
argument $x$ varies the absolute error function $\Delta (x)$ takes the
nonzero values $\lambda_1,-\lambda_2,\dots,(-1)^{m+n+1}\lambda_{m+n+2}$
having alternating
signs. In this case we shall say that {\it alternation} appears. If
$\vert\lambda_1\vert =\vert\lambda_2\vert =\dots =\vert\lambda_{m+n+2}\vert$,
then this alternation is {\it Chebyshev's}.
Denote by $\Delta_{\min}$ the best
possible absolute error of approximants of the form~(1) to the function
$f(x)$ (the numbers $m$ and $n$ are fixed). Suppose
$\lambda=\min\{\vert\lambda_1\vert ,\dots ,\vert\lambda_{m+n+2}\vert\}$. Then,
due to the generalized de la Vall\'ee--Poussin theorem, the inequality
$\Delta_{\min}\geq\lambda$ is valid; thus
$$
\Delta\geq\Delta_{\min}\geq\lambda.\tag11
$$

It is clear that $\Delta$ coincides with the greatest (in absolute value)
extremum of the function $f(x)$, and one can take the least
(in absolute value) extremum of this function for $\lambda$ (up to a sign).
The quantity
$$
q=\lambda /\Delta\tag12
$$
characterizes the proximity of the error of the given approximant
to the error of the
best approximant. It is clear that $0<q\leq 1$ and $q=1$ if the given
approximant is best. The closer the quantity $q$ to $1$ the higher the
approximant quality. From~(11) and~(12) it follows that
$$
\Delta_{\min}\geq q\cdot\Delta.\tag13
$$
Usually, the estimate~(13) is rather rough. The appearance of the
alternation itself indicates to the closeness of the error of the given
approximant to the best one, and the quantity $\Delta_{\min}/\Delta$
is, in general, much greater than the value of $q$.

Similarly, the quality of an approximant with respect to the best
relative error is evaluated.

If we can calculate the values of the approximated function for all the points
of the segment $[A,B]$ (or for a sufficiently ``dense'' set of
such points), and if the coefficients of the rational approximant $R(x)$ are
already known, then it is not hard to determine the points of local
extremum of the error function and to calculate the quantities
$\lambda_1,\lambda_2,\dots ,\lambda_{m+n+2}$, and also the quantities
$\lambda$ and $q$ by means of a special standard subroutine. The same
subroutine is also necessary for the construction of the best approximants by
means of an iteration method. A good program package for the
construction of rational approximants must contain a subroutine
of this sort as well
as a good subroutine for solving systems of linear algebraic equations and
must have, as a component part, routines which
implement both the algorithms for approximate solving the
approximation problem and the construction algorithms for best approximants.

\head\S5. Chebyshev polynomials and polynomial approximations\endhead

Chebyshev polynomials play an important role in approximation theory
and in computational practice (see, for example,~\cite{12, 13, 17, 18, 24, 33,
38}). We shall consider Chebyshev polynomials of the first kind.

These polynomials were defined by P.~L.~Chebyshev in the form
$$
T_n(x)=\cos(n\arccos x),\tag14
$$
where $n=0,1,\dots$. Assume that $\varphi =\arccos x$; representing
$\cos n\varphi$ via $\sin\varphi$ and $\cos\varphi$, it is not difficult
to verify that the right--hand side of formula~(14) coincides indeed
with a certain polynomial. In particular,
$$
\align
T_0(x)&=\cos 0=1,\\
T_1(x)&=\cos\varphi =\cos(\arccos x)=x,\\
T_2(x)&=\cos 2\varphi ={\cos}^2\varphi -{\sin}^2\varphi\\
&=T_1^2(x)-(1-T_1^2(x))=2x^2-1
\endalign
$$
and so on. For an actual computation of $T_n(x)$ the recurrence relation
$$
T_n(x)=2xT_{n-1}(x)-T_{n-2}(x)\tag15
$$
is usually used. Sometimes it is more convenient to consider the
polynomials $\overline T_n(x)=2^{-n+1}T_n(x)$ since the coefficient
at $x^n$ of the polynomial $\overline T_n(x)$ is equal to~$1$. The
polynomials mentioned above satisfy the recurrence relation
$$
\overline T_n(x)=x\overline T_{n-1}(x)-{1\over
4}\overline T_{n-2}(x).\tag16
$$
Consider a particular case of the problem of the best approximation,
the approximant to the function $f(x)=x^n$ on the segment $[-1,1]$
being looked for in the form of a polynomial $P(x)$ of degree $n-1$.
From the de la Vall\'ee--Poussin theorem it follows that the
approximant in question has the form $P(x)=x^n-\overline T_n(x)$.
In this case the
error function $\Delta(x)$ coincides with $\overline T_n(x)$ and
one can explicitly obtain the Chebyshev alternation points:
$x_k=-\cos{k\pi\over n}$, where $k=0,1,\dots ,n$. Indeed,
$$
\align
\overline T_n(x_k)&= 2^{-n+1}\cos n\big(\pi -{k\pi\over n}\big)\\
&=2^{-n+1}\cos (n-k)\pi ={(-1)^{n-k}\over 2^{n-1}},
\endalign
$$
i.e., $\overline T_n(x)$ takes its maximum value $1/2^{n-1}$ with
alternate signs at the points indicated above.
This implies an important consequence:
the best polynomial approximant of degree $n-1$ to the polynomial
$a_0+a_1x+\dots +a_nx^n$ on the segment $[-1,1]$ has the form
$a_0+a_1x+\dots +a_nx^n-a_n\overline T_n(x)$. This result allows to
reduce the degree of a polynomial (for example, of some polynomial
approximant) with a minimum loss of accuracy. The reduction of the
polynomial degree by means of successively applying the method indicated
above is called {\it economization}. The economization method is due to
C.~Lanczos, see~\cite{38}.

The monomials $x^0,x^1,\dots ,x^m$ can be expressed via the Chebyshev
polynomials $T_0,T_1,\dots ,T_m$. For $m>0$ the following formula is valid:
$$
x_m=2^{1-m}\sum_{k=0}^{[m/2]}a_k{m\choose k}T_{m-2k}(x),\tag17
$$
where ${m\choose k}$ are the binomial coefficients, $[m/2]$ is the integer
part of
the number $m/2$, $a_k=1/2$ for $k=m/2$ and $a_k=1$ for $k\neq m/2$.
The expansion of the polynomial $T_m$ in powers of $x$ for $m>0$ is
given by the formula
$$
T_m(x)={1\over 2}m\sum_{k=0}^{[m/2]}{(-1)^k(m-k-1)!\over
(m-2k)!}(2x)^{m-2k}.\tag18
$$
Finally, $x^0=T_0=1$. It is clear that the set of polynomials of the form
$\sum_{i=0}^nc_iT_i$, where $c_i$ are numerical coefficients, coincides with
the set of all polynomials of degree $n$.

The economization procedure mentioned above can also be described in
the following way. The initial polynomial $\sum_{i=0}^n a_ix^i$ can be
represented by means of formula~(17) in the form $\sum_{i=0}^nc_iT_i$. The
polynomial of degree $k$ obtained as the result of economization
coincides  with $\sum_{i=0}^kc_iT_i$. For functions represented in the
form of power series $f(x)=\sum_{i=0}^\infty a_ix^i$ it is not difficult
to obtain, by means of the economization method, polynomial approximants
on the segment $[-1,1]$ close to the best ones. For this purpose it is
necessary to replace $f(x)$ by its truncated Taylor series at the point
$x=0$, i.e., by the
polynomial $\sum_{i=0}^n a_ix^i$ approximating this function with a
high degree of accuracy, and then to obtain, by means of the economization
of this polynomial, the polynomial $\sum_{i=0}^kc_iT_i$ of the given degree
$k$. As $n\to\infty$, the quantity $\sum_{i=0}^kc_iT_i$ tends to the sum of
the first $k+1$ terms of the expansion of $f(x)$ into Fourier series with
respect to Chebyshev polynomials.

Denote by $L^2_w$ the Hilbert space of square integrable (with respect to the
measure $w(x)\,dx$) functions on the segment $[-1,1]$.
Suppose $w(x)=\sqrt{1-x^2}$. It
is not hard to verify that the Chebyshev polynomials $T_n$ form an
orthogonal (but not orthonormal) basis in $L_w^2$. The expansion
$$
f(x)=\sum_{i=0}^\infty c_iT_i\tag19
$$
of a function $f(x)$ into the series in Chebyshev polynomials (the
Fou\-rier--Chebyshev series) is easily reduced to the expansion of the
function
$f(\cos x)$ into the standard Fourier series in cosines. Among the
polynomials of degree $n$, the polynomial $P_n(x)=\sum_{i=0}^nc_iT_i$
gives the best approximation to the function $f(x)$ in $L_w^2$. The
following
result shows that this approximant on the segment $[-1,1]$ is close to the
best one in the sense of the absolute error.

\proclaim{Cheney theorem} Let $\varphi (x)$ be a function integrable on the
segment $[-1,1]$. If for $i=0,1,2,\dots ,k$ the equality
$$
\int\limits_{-1}^1\varphi (x)T_i(x)w(x)\,dx=0\tag20
$$
is valid, then $\varphi (x)$ either changes its sign in $[-1,1]$ at least
$k+1$ times or vanishes almost everywhere.\endproclaim

\demo{Proof} Assume that $\varphi (x)$ has exactly $m$ sign changes at the points
$x_1,\dots ,x_m$, where $0\leq m\leq k$. Suppose
$P(x)=\Pi_{i=1}^m(x-x_i)$;
since $P(x)$ is a polynomial of degree $m$,  we see that it can be
represented in the form of a linear combination
of the polynomials $T_0,T_1,\dots ,T_m$. Thus~(20) implies
$\int_{-1}^1\varphi (x)P(x)w(x)\,dx=0$. It is clear that the function
$\varphi (x)P(x)$ has no sign changes; thus the equality just obtained
means that
$\varphi (x)$ vanishes almost everywhere. The latter
proves the theorem.\enddemo

In a more general case, this result is proved in~\cite{39, p.110}. The proof
given above allows a generalization to the case of systems of
orthogonal polynomials of a sufficiently general form and arbitrary
segments of integration (including infinite ones), see~\cite{40}.

Now let us return to the function~(19) and to its approximant
$P_n(x)=\sum_{i=0}^nc_iT_i$. The absolute error function
$$
\Delta (x)=f(x)-P_n(x)=\sum_{i=n+1}^\infty c_iT_i
$$
is orthogonal to the polynomials $T_0,T_1,\dots ,T_n$, i.e.,
$$
\int_{-1}^1\Delta (x)T_i(x)w(x)\,dx=0
$$
for all $i=0,1,\dots ,n$. If the function $f(x)$ is continuous, then the
error function $\Delta (x)$ is also continuous. In this case from the Cheney
theorem it follows that either $\Delta (x)$ is identically zero
or has $n+1$ sign changes. This means that alternation is present,
i.e., the approximant $P_n(x)$ is close to the best one, and their proximity
can be evaluated by means of relations~(11)---(13).

While the truncated Taylor series $\sum_{i=0}^n{1\over n!}f^{(n)}(0)x^n$
gives the best approximant only in a neighborhood of the origin, the truncated
Fourier--Chebyshev series $\sum_{i=0}^nc_iT_i$ for the function $f(x)$ with
the same number of terms gives an approximant which is close to best one on
the entire segment $[-1,1]$.

The change of the variable $x\mapsto {1\over 2}[(B-A)x+A+B]$
reduces the problem of approximation on an arbitrary finite segment
$[A,B]$ to the case of the segment $[-1,1]$. Further, as a rule, we
shall consider the latter case.

\head\S6. Ill-conditioned problems and rational approximations\endhead

Let $\{\varphi_0,\varphi_1,\dots ,\varphi_n\}$ and
$\{\psi_0,\psi_1,\dots ,\psi_m\}$ be collections consisting of linearly
independent functions of the argument $x$ belonging some (possibly
multidimensional) set $X$. Consider the problem of constructing
an approximant of the form
$$
R(x)={a_0\varphi_0+a_1\varphi_1+\dots +a_n\varphi_n\over
b_0\psi_0+b_1\psi_1+\dots +b_m\psi_m}\tag21
$$
to a given function $f(x)$ defined on $X$. If $X$ coincides with a real line
segment $[A,B]$, $\varphi_k=x^k$ and $\psi_k=x^k$ for
all $k$, then the expression~(21) turns out to be a rational function
of the form~(1)
(see the Introduction). It is clear that expression~(21) also gives a
rational function in the case when we take Chebyshev polynomials
$T_k$ or, for example, Legendre, Laguerre, Hermite, etc.
polynomials as $\varphi_k$ and $\psi_k$.

Fix an abstract construction method for an approximant of the
form (21)
and consider the problem of computing the coefficients $a_i$, $b_j$.
Quite often this problem is ill-conditioned, i.e., small perturbations of the
approximated function $f(x)$ or a calculation errors lead to considerable
errors in the values of coefficients. For example, the problem of computing
coefficients for best rational approximants (including polynomial
approximants) for high degrees of the numerator or the denominator
is ill-conditioned.

The instability with respect to the calculation error can be related
both to the abstract construction method of approximation (i.e.,
with the formulation of the problem) and to the particular algorithm
implementing the method. The fact that the problem of computing
coefficients for the best approximant is ill-conditioned is related
to the formulation of this problem.
This is also valid for other construction methods for rational
approximants with a sufficiently large number of coefficients.
But an unfortunate choice of the
algorithm implementing a certain method can aggravate troubles
connected with ill-conditioning.

Several construction methods for approximants of the form~(21)
are connected with solving systems of linear algebraic equations. This
procedure can lead to a large error if the corresponding matrix is
ill-conditioned. Consider an arbitrary system of linear algebraic equations
$$
Ay=h,\tag22
$$
where $A$ is a given square matrix of order $N$ with components
$a_{ij}$ $(i,j=1,\dots ,N)$, $h$ is a given vector column with components
$h_i$, and $y$ is an unknown vector column with
components $y_i$. Define the vector norm by the equality
$$
\Vert y\Vert =\sum_{i=1}^N\vert x_i\vert\tag23
$$
(this norm is more convenient for calculations than
$\sqrt{x_1^2+\dots +x_N^2}$). Then the matrix norm is determined by the
equality
$$
\Vert A\Vert =\max_{\Vert y\Vert =1}\Vert Ay\Vert =\max_{1\leq j\leq N}
\sum_{i=1}^N\vert a_{ij}\vert.\tag24
$$

If a matrix $A$ is nonsingular, then the quantity
$$
\cond (A)=\Vert A\Vert\cdot\Vert A^{-1}\Vert\tag25
$$
is called the {\it condition number} of the matrix $A$ (see, for example,
\cite{41}). Since $y=A^{-1}h$, we see that the absolute error $\Delta y$
of the vector $y$ is connected with the absolute error of the vector $h$
by the relation $\Delta y=A^{-1}\Delta h$, whence
$$
\align
\Vert\Delta y\Vert&\leq\Vert A^{-1}\Vert\cdot\Vert\Delta h\Vert\\
\text{ and }
\Vert\Delta y\Vert /\Vert y\Vert&\leq\Vert A^{-1}\Vert\cdot
(\Vert h\Vert /\Vert y\Vert)(\Vert\Delta h\Vert /\Vert h\Vert).
\endalign
$$
Taking into account the fact that
$\Vert h\Vert\leq\Vert A\Vert\cdot\Vert y\Vert$, we finally obtain
$$
\Vert\Delta y\Vert /\Vert y\Vert\leq\Vert A\Vert\cdot\Vert A^{-1}
\Vert\cdot\Vert\Delta h\Vert /\Vert h\Vert ,\tag26
$$
i.e., the relative error of the solution $y$ is estimated via the relative
error of the vector $h$ by means of the condition number. It is clear
that~(26) can turn into an equality. Thus, if the condition number is of
order $10^k$, then, because of round--off errors in $h$, we
can lose $k$ decimal digits of $y$.

Similarly, the contribution of the error of the matrix $A$
is evaluated. Finally, the dependence of $\cond (A)$ on the choice of a
norm is weak.
A method of rapid estimation of the condition number is described
in~\cite{41, \S 3.2}. The analysis of the cases when the condition number
gives a much too pessimistic error estimate is given in~\cite{42}.

As an example, we note that the coefficients of the polynomial $P_n(x)$ which
give the best approximant to the function $f(x)$ in the metric of the
Hilbert
space $L_w^2$ (see~\S5 above) can be determined from the system of
equations
$$
\int\limits_{-1}^1\big(f(x)-P_n(x)\big)x^kw(x)\,dx=0,\tag27
$$
where $k=0,1,\dots ,n$. With respect to coefficients of the polynomial
$P_n(x)$ (in powers of $x$ or in Chebyshev polynomials) these equations
are linear and algebraic. But due to the fact that the
monomials $x^k$ are ``almost linearly dependent'', system~(27) is very
ill-conditioned. The equi\-valent system
$$
\int\limits_{-1}^1\big(f(x)-P_n(x)\big)T_k(x)w(x)\,dx=0\tag"$(27')$"
$$
is better conditioned, but in this case it is also preferable to use the
economization procedure or to determine the coefficients $c_i$ in~(19)
by formulas
$$
\aligned
c_0&={1\over\pi}\int\limits_{-1}^1f(x)w(x)\,dx,\\
c_i&={2\over\pi}\int\limits_{-1}^1f(x)T_i(x)w(x)\,dx,\qquad i\geq 1.
\endaligned
\tag28
$$
We recall that here $w(x)=1/\sqrt{1-x^2}$.

\head\S7. The effect of error autocorrection\endhead

Fix an abstract construction method (problem) for an approximant of the
form (21) to the function $f(x)$. Let the coefficients $a_i$, $b_j$
give an exact or an approximate solution of this problem, and
let the $\tilde a_i$, $\tilde b_j$ give another approximate solution
obtained in the same way. Denote by $\Delta a_i$, $\Delta b_j$ the
absolute errors of the coefficients, i.e., $\Delta a_i=\tilde a_i-a_i$,
$\Delta b_j=\tilde b_j-b_j$; these errors arise due to perturbations
of the approximated function $f(x)$ or due to calculation errors.
Set
$$
\alignedat2
&P(x)=\sum_{i=0}^na_i\varphi_i,&&\qquad Q(x)=\sum_{j=0}^mb_j\psi_j,\\
&\Delta P(x)=\sum_{i=0}^n\Delta a_i\varphi_i,&&\qquad\Delta
Q(x)=\sum_{j=0}^m\Delta b_j\psi_j,\\
&\widetilde{P}(x)=P+\Delta P,&&\qquad\widetilde{Q}(x)=Q+\Delta Q.
\endalignedat
$$

It is easy to verify that the following exact equality is valid:
$$
{P+\Delta P\over Q+\Delta Q}-{P\over Q}={\Delta Q\over Q}\Bigg({\Delta
P\over
\Delta Q}-{P\over Q}\Bigg).\tag29
$$

As mentioned in the Introduction, the fact that the problem of calculating
coefficients is ill-conditioned can nevertheless be accompanied by
high accuracy of the
approximants obtained. This means that the approximants $P/Q$ and
$\widetilde{P}/\widetilde{Q}$ are close to the approximated function and,
therefore, are close to each other, although the coefficients of these
approximants differ greatly. In this
case the relation $\Delta Q/\widetilde{Q}=\Delta Q/(Q+\Delta Q)$ of
the denominator considerably exceeds in absolute value the left-hand
side of equality~(29). This is possible only in the case when the
difference $\Delta P/\Delta Q-P/Q$ is small, i.e., the function
$\Delta P/\Delta Q$ is close to $P/Q$, and, hence, to the approximated
function. Thus the function $\Delta P/\Delta Q$ will be called the
{\it error approximant}.
For a special case, this concept was actually introduced in~\cite{5}.
In the sequel, we shall see that in many cases the error approximant
provides indeed a good
approximation for the approximated function, and, thus, $P/Q$ and
$\widetilde{P}/\widetilde{Q}$ differ from each other by a product of
small quantities in the right-hand side of~(29). The thing is that
the errors $\Delta a_i$, $\Delta b_j$ are not arbitrary, but are connected by
certain relations.

Let an abstract construction method for the approximant of the
form (21)
be linear in the sense that the coefficients of the approximant can be
determined from a homogeneous system of linear algebraic equations. The
homogeneity condition is connected with the fact that, when multiplying
the numerator and the denominator of fraction~(21) by the same nonzero
number, the approximant~(21) does not change. Denote by $y$ the vector whose
components are the coefficients $a_0, a_1,\dots ,a_n$,
$b_0,b_1,\dots ,b_m$.
Assume that the coefficients can be obtained from the homogeneous system
of equations
$$
Hy=0,\tag30
$$
where $H$ is a matrix of dimension $(m+n+2)\times (m+n+1)$.

The vector $\tilde y$ is an approximate solution of system~(30) if the
quantity $\Vert H\tilde y\Vert$ is small. If $y$ and $\tilde y$ are
approximate solutions of system~(30), then the vector $\Delta y=\tilde y-y$
is also an approximate solution of this system since
$\Vert H\Delta y\Vert=\Vert H\tilde y-Hy\Vert\leq\Vert H\tilde y\Vert
+\Vert Hy\Vert$. Thus it is natural to assume that the function
$\Delta P/\Delta Q$ corresponding to the solution $\Delta y$ is
an approximant to $f(x)$. It is clear that the order of the residual
of the approximate solution $\Delta y$ of system~(30), i.e., of the
quantity $\Vert H\Delta y\Vert$, coincides with the order of the largest
of the residuals of the approximate solutions $y$ and $\tilde y$. For a fixed
order of the residual the increase of
the error $\Delta y$ is compensated by the fact that $\Delta y$
satisfies the system of equations~(30)
with greater ``relative'' accuracy, and
the latter, generally speaking, leads to the increase of
accuracy of the error approximant.

To obtain a certain solution of system~(30), one usually adds to this system
a normalization condition of the form
$$
\sum_{i=0}^n\lambda_ia_i+\sum_{j=0}^m\mu_jb_j=1,
\tag31
$$
where $\lambda_i$, $\mu_j$ are numerical coefficients. As a rule, the
equality $b_0=1$ is taken as the normalization condition (but this
is not always successful with respect to minimizing the calculation
errors).

Adding equation~(31) to system~(30), we obtain a nonhomogeneous system of
$m+n+2$ linear algebraic equations of type~(22). If the approximate
solutions $y$ and $\tilde y$ of system~(30) satisfy condition~(31), then
the vector $\Delta y$ satisfies the condition
$$
\sum_{i=0}^n\lambda_i\Delta a_i+\sum_{j=0}^m\mu_j\Delta b_j=0,
\tag"$(31')$"
$$

It is clear that the above reasoning is not rigorous; for each
specific construction method for approximations it is necessary to carry
out some additional analysis. More accurate reasoning is given below, in~\S8,
for the classical Pad\'e approximants, and in~\S14, for the linear and
nonlinear
Pad\'e--Chebyshev approximants. The presence of the error autocorrection
mechanism described above is also verified by a numerical experiment
(see below).

The effect of error autocorrection reveals itself for certain nonlinear
construction methods for rational approximations as well.
One of these methods is
considered below, in \S12--14 (nonlinear Pad\'e--Chebyshev approximation).

It must be emphasized that (as noted in \S3) the coefficients
of the best Chebyshev approximant satisfy the system of linear algebraic
equations~(9) and are computed as approximate solutions of this system on
the last step of the iteration process in algorithms of Remez's type.
Thus, the construction methods for the best rational approximants can be
regarded as linear. At least for some functions (say, for
$\cos\pi /4x$, $-1\leq x\leq 1$) the linear and the nonlinear
Pad\'e--Chebyshev
approximants are very close to the best ones in the sense of the relative
and the absolute errors, respectively. The results that arise when applying
calculation algorithms for Pad\'e--Chebyshev approximants can be regarded as
approximate solutions of system~(9) which determines the best approximants.
Thus the presence of the effect of error autocorrection for Pad\'e--Chebyshev
approximants gives an additional argument in favor of the conjecture
that this effect also
takes place for the best approximants.

Finally, note that the basic relation~(29) becomes meaningless if one seeks
an approximant in the form
$a_0\varphi_0+a_1\varphi_1+\dots +a_n\varphi_n$, i.e., the denominator
in~(21)
is reduced to $1$. However, in this case the effect of error
autocorrection (although much weakened) is also possible; this is
connected with
the fact that the errors $\Delta a_i$ approximately satisfy certain
relations. Such a situation can arise when using the least
squares method.

\head\S8. Pad\'e approximations\endhead

Let the expansion of a function $f(x)$ into a power series (the Taylor
series at zero) be given, i.e.,
$$
f(x)=\sum_{i=0}^\infty c_ix^i.\tag32
$$
The classical Pad\'e approximant for $f(x)$ is a rational function of the form
$$
R(x)=P_n(x)/Q_m(x),\tag33
$$
where $P_n(x)$ and $Q_m(x)$ are polynomials of degree $n$ and $m$,
respectively, satisfying the relation
$$
Q_m(x)f(x)-P_n(x)=O(x^{m+n+1}).\tag34
$$
Let
$$
\aligned
P_n(x)&=a_0+a_1x+\dots +a_nx^n,\\
Q_m(x)&=b_0+b_1x+\dots +b_mx^m.
\endaligned
\tag35
$$
If $b_0\neq 0$, then~(34) means that
$$
f(x)-P_n(x)/Q_m(x)=O(x^{m+n+1}),\tag"$(34')$"
$$
i.e., the first $m+n+1$ terms of the Taylor expansion in powers of
$x$ (to $x^{m+n}$ inclusive) of $f(x)$ and $R(x)$ are the same. The
Pad\'e approximation gives the best approximant in a small neighborhood
of zero; it is a natural generalization of the expansion of functions into
Taylor series and is closely connected with the expansion of functions into
continued fractions. Numerous papers are devoted to the Pad\'e
approximation; see, for example,~\cite{11--16, 5, 6}.

One can evaluate the coefficients $b_j$ in the denominator of fraction~(33)
by solving the homogeneous system of linear equations
$$
\sum_{j=1}^mc_{n+k+j}b_j=-b_0c_{n+k},\tag36
$$
where $k=1,\dots ,m$ and $c_l=0$ for $l<0$. One can take any nonzero
constant as $b_0$. The coefficients $a_i$ are given by the
formulas
$$
a_i=\sum_{k=0}^ib_kc_{i-k}=\sum_{k=0}^ib_{i-k}c_k.\tag37
$$
The text of the corresponding Fortran program is given in~\cite{11}.

For large $m$ the system~(36) is ill-conditioned. Moreover, the problem of
computation for coefficients of Pad\'e approximants is also
ill-conditioned independent of a particular solving algorithm for this
problem, see~\cite{6, 43, 44}. In Y.~L.~Luke's paper~\cite{5} the
following reasoning is given. Let $\Delta a_i$, $\Delta b_j$ be the
errors in the coefficients $a_i$, $b_j$ which arise when numerically
solving system~(36).
We shall ignore the errors of the quantities $c_i$ and $x$ and
we shall consider that, according to~(37), the errors in the
coefficients $a_i$ have the form
$$
\Delta a_i=\sum_{k=0}^i\Delta b_{i-k}c_k.\tag"{($37'$)}"
$$
From~(37') it follows that
$$
\align
f\Delta Q-\Delta P&=\sum_{j=0}^m\Delta b_jx^j\sum_{i=0}^\infty c_ix^i-
\sum_{i=0}^n\Delta a_ix^i\\
&=\sum_{i=0}^m\sum_{i=0}^\infty\Delta b_jc_ix^{i+j}-\sum_{i=0}^n
\sum_{k=0}^i\Delta b_{i-k}c_kx^i,
\endalign
$$
the latter, after the change of indices, yields the relation
similar to~(34):
$$
\Delta Q(x)f(x)-\Delta P(x)=O(x^{n+1}).
$$
Thus, there are reasons to expect
that the error approximant approximates indeed the function $f(x)$ and the
effect of error autocorrection takes place. In~\cite{5} the corresponding
experimental data for the function $e^{-x}$ for $x=2$, $m=n=6,7,\dots , 14$,
and for $x=5$ are given and the experiments with the functions
$x^{-1}\ln (1+x)$, $(1+x)^{\pm 1/2}$, $xe^x\int_x^\infty t^{-1}e^{-t}dt$
are briefly described; see also~\cite{6}.

A natural generalization of the classical Pad\'e approximant is the
multipoint
Pad\'e approximant (or Pad\'e approximant of the second kind), i.e., a
rational function of the form~(33) whose values
coincide with values of the approximated function $f(x)$
at some points $x_i$ ($i=1,2,\dots ,m+n+1$).
This definition is extended to the case of multiple points,
and for $x_i=0$ for all $i$ it leads to the classical Pad\'e
approximations
see~\cite{11, 14, 15}. The calculation of coefficients in the multipoint
Pad\'e approximant can be reduced to solving a system of linear
equations, and there are reasons to suppose that in this case the
effect of error autocorrection takes place as well.

\head\S9. Linear Pad\'e--Chebyshev approximations and the
PADE program\endhead

Consider the approximant of the form~(33) to the function $f(x)$ on the
segment $[-1,1]$. The absolute error function of this approximant has the
following form:
$$
\Delta (x)=\Phi (x)/Q_m(x),
$$
where
$$
\Phi (x)=f(x)Q_m(x)-P_n(x).\tag38
$$
The function $R_{m,n}(x)=P_n(x)/Q_m(x)$ is called the {\it linear
Pad\'e--Chebyshev approximant to the function} $f(x)$ if
$$
\int\limits_{-1}^1\Phi (x)T_k(x)w(x)\,dx=0,\qquad k=0,1,\dots ,m+n,\tag39
$$
where $T_k(x)$ are the Chebyshev polynomials,
$w(x)=1/\sqrt{1-x^2}$. This concept (in a different form) was introduced
in~\cite{45} and allows a generalization to the case of other orthogonal
polynomials (see~\cite{11, 33, 34, 39, 40}).  Approximants of this kind always
exist~\cite{39}. Reasoning in the same way as in~\S5 and applying
Cheney's theorem, we can find out why the linear
Pad\'e--Chebyshev approximants are close to the best ones.

Let $P_n(x)$ and $Q_m(y)$ be represented in the form~(35). Then the system of
equations~(39) is equivalent to the following system of linear algebraic
equations with respect to the coefficients $a_i$, $b_j$:
$$
\sum_{j=0}^mb_j\int\limits_{-1}^1{x^jT_k(x)f(x)\over\sqrt{1-x^2}}\,dx
-\sum_{i=0}^na_i\int\limits_{-1}^1{x^iT_k(x)\over\sqrt{1-x^2}}\,dx=0.\tag40
$$

The homogeneous system~(40) can be transformed into a nonhomogeneous one by
adding a normalization condition; in particular, any of the following
equalities can be taken as this condition:
$$
\alignat 3
&b_0&=1, \tag41\\
&b_m&=1, \tag42\\
&a_m&=1. \tag43
\endalignat
$$

In~\cite{1, 9} the program PADE (in Fortran, with double precision) which
allows to construct rational approximants by solving the system of
equations of type~(40) is briefly described.
The complete text of a certain version of this program and its detailed
description can be found in the Collection of algorithms
and programs of the Research Computer Center of the Russian Acad.
Sci~\cite{7}. For even functions the approximant is looked for in the form
$$
R(x)={a_0+a_1x^2+\dots +a_n(x^2)^n\over b_0+b_1x^2+\dots
+b_m(x^2)^m},\tag44
$$
and for odd functions it is looked for in the form
$$
R(x)=x{a_0+a_1x^2+\dots +a_n(x^2)^n\over b_0+b_1x^2+\dots
+b_m(x^2)^m},\tag45
$$
respectively.
The program computes the values of coefficients of the approximant, the
absolute and the relative errors, and gives the information which
allows to estimate
the quality of the approximation (see~\S4 above). In particular, a version
of the PADE program is implemented by means of minicomputer of
SM--4 class constructs the error curve,
determines the presence of alternation, and produces the estimate of the
quality of the approximation by means of quantity~(12). Using a
subroutine, the user introduces the function defined by means of any
algorithm on an
arbitrary segment~$[A,B]$, introduces the boundary points of this segment, the
numbers $m$ and $n$, and the number of control parameters. In particular,
one can choose the normalization condition of type~(41)--(43), look for an
approximant in the form~(44) or~(45) and so on. The change of the
variable reduces the approximation on any segment $[A,B]$ to the
approximation on
the segment $[-1,1]$. Therefore, we shall consider the case
when $A=-1$, $B=1$ in the sequel unless otherwise stated.

For the calculation of integrals, the Gauss--Hermite--Chebyshev quadrature
formula is used:
$$
\int\limits_{-1}^1{\varphi (x)\over\sqrt{1-x^2}}\,dx={\pi\over s}
\sum_{i=1}^s\varphi\Big(\cos{2i-1\over 2s}\pi\Big),\tag46
$$
where $s$ is the number of interpolation points;
for polynomials of degree $2s-1$
this formula is exact, so that the precision of formula~(46)
increases rapidly as the parameter $s$ increases and depends on the
quality of the approximation of the function $\varphi (s)$ by polynomials.
To calculate the values of Chebyshev polynomials, recurrence relation~(15)
is applied.

If the function $f(x)$ is even and an approximant is looked for the
form~(44),
then system~(40) is transformed into the following system of equations:
$$
\sum_{i=0}^na_i\int\limits_{-1}^1{x^{2i}T_{2k}(x)\over\sqrt{1-x^2}}\,dx
-\sum_{j=
0}^mb_j\int\limits_{-1}^1{x^{2j}T_{2k}(x)f(x)\over\sqrt{1-x^2}}\,dx=0,
\tag47
$$
where $k=0,1,\dots ,m+n$. If $f(x)$ is an odd function and an approximant
is looked for in the form~(45), then, first, by means of the solution of
system~(47) complemented by one of the normalization conditions, one
determines an approximant of the form~(44) to the even
function $f(x)/x$, and then the
obtained approximant is multiplied by $x$. This procedure allows
to avoid a large relative error for $x=0$.

The possibilities of the PADE program are demonstrated in Table~1. This
table contains errors of certain approximants obtained by means of this
program. For every approximant, the absolute error $\Delta$, the relative
error $\delta$, and (for comparison) the best possible relative
error $\delta_{\min}$ taken from~\cite{3} are indicated. The function
$\sqrt{x}$ is approximated on the segment $[1/2,1]$ by the expression of
the
form~(1), the function $\cos{\pi\over 4}x$ is approximated on the segment
$[-1,1]$ by the expression of the form~(44), and all the others are
approximated on the same segment by the expression of the form~(45).

\centerline{Table 1}

$$
\vmatrix
\format \c&\vrule\quad\c\quad&\vrule\quad\c\quad&\vrule\quad\c\quad&
\vrule\quad\c\quad&\vrule\quad\c\\
\text{ Function }&m&n&\Delta&\delta&\delta_{\min}\cr
\_\_\_\_\_\_\_\_\_\_\_&\_\_\_&\_\_\_&\_\_\_\_\_\_\_\_\_\_\_\_&
\_\_\_\_\_\_\_\_\_\_\_\_&\_\_\_\_\_\_\_\_\_\_\_\_\cr
\sqrt{x}&2&2&0.8\cdot 10^{-6}&1.13\cdot 10^{-6}&0.6\cdot 10^{-6}\cr
\sqrt{x}&3&3&1.9\cdot 10^{-9}&2.7\cdot 10^{-9}&1.12\cdot 10^{-9}\cr
\cos{\pi\over 4}x&0&3&0.28\cdot 10^{-7}&0.39\cdot 10^{-7}&0.32\cdot
10^{-7}\cr
\cos{\pi\over 4}x&1&2&0.24\cdot 10^{-7}&0.34\cdot 10^{-7}&0.29\cdot
10^{-7}\cr
\cos{\pi\over 4}x&2&2&0.69\cdot 10^{-10}&0.94\cdot 10^{-10}&0.79\cdot
10^{-10}\cr
\cos{\pi\over 4}x&0&5&0.57\cdot 10^{-13}&0.79\cdot 10^{-13}&0.66\cdot
10^{-13}\cr
\cos{\pi\over 4}x&2&3&0.4\cdot 10^{-13}&0.55\cdot 10^{-13}&0.46\cdot
10^{-13}\cr
\sin{\pi\over 4}x&0&4&0.34\cdot 10^{-11}&0.48\cdot 10^{-11}&0.47\cdot
10^{-11}\cr
\sin{\pi\over 4}x&2&2&0.32\cdot 10^{-11}&0.45\cdot 10^{-11}&0.44\cdot
10^{-11}\cr
\sin{\pi\over 4}x&0&5&0.36\cdot 10^{-14}&0.55\cdot 10^{-14}&0.45\cdot
10^{-14}\cr
\sin{\pi\over 2}x&1&1&0.14\cdot 10^{-3}&0.14\cdot 10^{-3}&0.12\cdot
10^{-3}\cr
\sin{\pi\over 2}x&0&4&0.67\cdot 10^{-8}&0.67\cdot 10^{-8}&0.54\cdot
10^{-8}\cr
\sin{\pi\over 2}x&2&2&0.63\cdot 10^{-8}&0.63\cdot 10^{-8}&0.53\cdot
10^{-8}\cr
\sin{\pi\over 2}x&3&3&0.63\cdot 10^{-13}&0.63\cdot 10^{-13}&0.5\cdot
10^{-13}\cr
\tg{\pi\over 4}x&1&1&0.64\cdot 10^{-5}&0.64\cdot 10^{-5}&0.57\cdot
10^{-5}\cr
\tg{\pi\over 4}x&2&1&0.16\cdot 10^{-7}&0.16\cdot 10^{-7}&0.14\cdot
10^{-7}\cr
\tg{\pi\over 4}x&2&2&0.25\cdot 10^{-10}&0.25\cdot 10^{-10}&0.22\cdot
10^{-10}\cr
\arctg x&0&7&0.75\cdot 10^{-7}&10^{-7}&10^{-7}\cr
\arctg x&2&3&0.16\cdot 10^{-7}&0.51\cdot 10^{-7}&0.27\cdot 10^{-7}\cr
\arctg x&0&9&0.15\cdot 10^{-8}&0.28\cdot 10^{-8}&0.23\cdot 10^{-8}\cr
\arctg x&3&3&0.54\cdot 10^{-9}&1.9\cdot 10^{-9}&0.87\cdot 10^{-9}\cr
\arctg x&4&4&0.12\cdot 10^{-11}&0.48\cdot 10^{-11}&0.17\cdot 10^{-11}\cr
\arctg x&5&4&0.75\cdot 10^{-13}&3.7\cdot 10^{-13}&0.71\cdot 10^{-13}\cr
\endvmatrix
$$

The PADE program is comparatively simple and compact; it includes the standard
subroutine DGELG for solving systems of linear algebraic equations (this
subroutine is taken from~\cite{46}) and a subroutine of numerical
integration, and also a number of service, test and auxiliary modules.
No additional software is used. The program needs minimum
hardware requirements and can be implemented by means of
any computer
having a Fortran compiler, random-access memory of sufficient
volume, a printer or a display.

The version of the PADE program described in~\cite{7} is implemented by means
of computers of IBM 360/370 class and requires 60 K bytes of main memory;
the volume of this program in Fortran (including comments) is 581 lines
(cards).
The program execution time depends on the type of the computer,
on the approximated function, and on the values of control parameters.
For example, the CPU time for determining, by means of the PADE program,
an approximant of the form~(1) to the function $\sqrt{x}$ on
the segment $[1/2,x]$ for $m=n=2$ is 4\.4s. In this case the
normalization~(43)
is applied, and the number of checkpoints used while
estimating the error is 1200; the compilation time
is not taken into account
\footnote{One can sufficiently decrease the number of
checkpoint without considerable loss of accuracy of error
estimation (in the present case, for example, to 200 points).}.

One of the versions of the program gives the estimate of the quality of
the  approximant obtained according to formula~(12) (see~\S4 above).
For example,
for the function $\sin{\pi\over 2}x$ for $m=n=2$, and for the relative
error we have $q=0.0625$ whence it follows that
$\delta_{\min}\geq q\delta\approx 0.4\cdot 10^{-9}$.
This estimate is rough and in fact, as is shown in Table~1,
$\delta_{\min}/\delta\approx 0.84$. For the
absolute error the program gives in this case $q=0.71$. The latter
indicates to the closeness of this error to the best possible. The version
of the program mentioned above allows to carry out the calculations in
interactive
mode varying the degrees $m$ and $n$, the boundary points of the segment
$[A,B]$, the branches of the algorithm, the number of checkpoints
when the errors are calculated,
the number of interpolation points in the quadrature formula~(46),
and to estimate rapidly the quality of the approximation according
to the error curve.

{\bf Remark.} The program of constructing classical Pad\'e approximants
gi\-ven in~\cite{11} is also called PADE, but, of course,
here and in~\cite{11} different programs are discussed.

\head\S10. The PADE program. Analysis of the algorithm\endhead

The quality of an approximant obtained by means of the PADE program
mainly depends on the behavior of the denominator of this approximant
and on the calculation errors. The fact that the corresponding systems of
algebraic equations are ill-conditioned is the most unpleasant the
source of errors of the method under consideration.
Seemingly, the methods of this kind are not widely used due
to this reason.

The condition numbers of systems of equations that arise while calculating,
by means of the PADE program, the approximants considered above are
also very large, for example, while calculating the approximant of the
form~(5) on the segment $[-1,1]$ to $\sin{\pi\over 2}x$ for $m=n=3$, the
corresponding condition number is of order $10^{13}$.
As a result, the coefficients of the approximant are
determined with a large error. In particular, a small perturbation of
the system of linear equations arising when passing from computer
ICL 4--50 to ES--1045 (because of the calculation errors)
gives rise to large perturbations in the coefficients of the approximant.
Fortunately, the effect of error autocorrection (see~\S7 above)
improves the situation, and the errors of the approximant have no
substantial changes under this perturbation. This fact is
described in the Introduction, where concrete examples are also given.

Consider some more examples connected with passing from
ICL 4--50 to ES--1045. The branch of the algorithm which corresponds
to the normalization condition~(41) (i.e., to $b_0=1$) is considered.
For $\arctg x$ the calculation of an approximant of the form~(45)
on the segment $[-1,1]$ for $m=n=5$ by means of ICL--4--50 computer
gives an approximation with the absolute error $\Delta=0.35\cdot 10^{-12}$
and the relative error $\delta =0.16\cdot 10^{-11}$. The corresponding
system of linear algebraic equations has the condition number of
order $10^{30}$! Passing to ES--1045 we obtain the
following: $\Delta =0.5\cdot 10^{-14}$, $\delta =0.16\cdot 10^{-12}$,
the condition number is of order $10^{14}$, and the errors
$\Delta a_1$ and $\Delta b_1$ in the coefficients $a_1$ and $b_1$ in~(45)
are greater in absolute value than $1$! This example shows that
the problem of computing condition number of an ill-conditioned system is,
in its turn, ill-conditioned. Indeed, the
condition number is, roughly speaking, determined by values of
coefficients of the inverse matrix (see~\S6 above, eqs~(24)
and~(25)), every column of the inverse matrix being the solution of the
system of equations with the initial matrix of coefficients, i.e., of
an ill-conditioned system.

Consider in more detail the effect of error autocorrection for the
approximant of the form~(44) on the segment $[-1,1]$ to the
function $\cos{\pi\over 4}x$ for $m=2$, $n=3$. Constructing this
approximant both on the ICL--4--50 and the ES--1045 computer
results in the approximation with the absolute error
$\Delta =0.4\cdot 10^{-13}$
and the relative error $\delta =0.55\cdot 10^{-13}$ which are close to
the best possible. In both the cases the condition number is of order $10^9$.
The coefficients of the approximants obtained by means of the
computers mentioned above
and the coefficients of the error approximant (see~\S7 above) are
as follows:
$$
\aligned
\tilde a_0&=0.9999999999999600,\\
a_0&=0.9999999999999610,\\
\tilde a_1&=-0.2925310453579570,\\
a_1&=-0.2925311264716216,\\
\tilde a_2&=10^{-1}\cdot 0.1105254254716866,\\
a_2&=10^{-1}\cdot 0.1105256585556549,\\
\tilde a_3&=10^{-3}\cdot 0.1049474500904401,\\
a_3&=10^{-3}\cdot 0.1049482094850086,\\
b_0&=1,\\
\tilde b_0&=1,\\
\tilde b_1&=10^{-1}\cdot 0.1589409217324021,\\
b_1&=10^{-1}\cdot 0.1589401105960337,\\
\tilde b_2&=10^{-3}\cdot 0.1003359011092697,\\
b_2&=10^{-3}\cdot 0.1003341918083529,
\endaligned
\qquad
\aligned
\Delta a_0&=-10^{-15},\\
{}\\
\Delta a_1&=10^{-7}\cdot 0.811136646,\\
{}\\
\Delta a_2&=-10^{-7}\cdot 0.2330839683,\\
{}\\
\Delta a_3&=10^{-9}\cdot 0.7593947685,\\
{}\\
\Delta b_0&=0,\\
{}\\
\Delta b_1&=10^{-7}\cdot 0.8111363684,\\
{}\\
\Delta b_2&=10^{-8}\cdot 0.17093009168.
\endaligned
$$

Thus, the error approximant has the form
$$
{\Delta P\over\Delta Q}={\Delta a_0+\Delta a_1x^2+\Delta a_2x^4+\Delta
a_3x^6
\over\Delta b_1x^2+\Delta b_2x^4}.\tag48
$$
If the relatively small quantity $\Delta a_0=-10^{-15}$
in~(48) is omitted,
then, as testing by means of a computer shows (2000 checkpoints), this
expression is an approximant to the function $\cos{\pi\over 4}x$ on the
segment $[-1,1]$ with the absolute and the relative errors
$\Delta =\delta =0.22\cdot 10^{-6}$.

But the polynomial $\Delta Q$ is zero at $x=0$, and the polynomial
$\Delta P$ takes a small, but nonzero value at $x=0$.
Fortunately, equality~(29) can be rewritten in the following way:
$$
{\widetilde{P}\over\widetilde{Q}}-{P\over Q}={\Delta P\over\widetilde{Q}}-
{\Delta Q\over\widetilde{Q}}\cdot{P\over Q}.\tag49
$$
Thus, as $\Delta Q\to 0$, the effect of error autocorrection
arises because the quantity $\Delta P$ is close to
zero, and the error of the approximant $P/Q$ is determined by the
error of the coefficient $a_0$. The same situation also take place when
the polynomial $\Delta Q$ vanishes at an arbitrary point
$x_0$ belonging to the segment $[A,B]$ where the function is
approximated. It is clear that if one chooses the standard normalization
($b_0=1$), then the error approximant has actually two coefficients less than
the initial one. Relations~(38) and~(39) show that in the
general case the normalization conditions $a_n=1$ or $b_m=1$ result
in the following: the coefficients of the error approximant form an
approximate solution of the homogeneous system of linear algebraic
equations whose exact solution determines the Pad\'e--Chebyshev approximant
having one coefficient less than the initial one. The effect of error
autocorrection improves again the accuracy of this error
approximant; thus,
``the snake bites its own tail''. A situation also
arises in the case when the approximant of the form~(44) to an
even function is constructed by solving the system of equations~(47).

Sometimes it is possible to decrease the error of the approximant by means of
the fortunate choice of the normalization condition. As
an example, consider the approximation of the function $e^x$
on the segment $[-1,1]$ by rational
functions of the form~(1) for $m=15$, $n=0$.
For the traditionally accepted normalization $b_0=1$, the PADE
program yields an
approximant with the absolute error $\Delta =1.4\cdot 10^{-14}$ and
the relative error $\delta =0.53\cdot 10^{-14}$. After passing to the
normalization condition $b_{15}=1$, the errors are reduced nearly one half:
$\Delta =0.73\cdot 10^{-14}$, $\delta =0.27\cdot 10^{-14}$.
Note that the condition number increases: in the first
case it is $2\cdot 10^6$, and in the second case it is $0\cdot 10^{16}$.
Thus the error decreases notwithstanding the fact that the system of
equations becomes drastically ill-conditioned. This example shows that
the increase of accuracy of the error approximant can be
accompanied by the increase
of the condition number, and, as experiments show, by the increase of errors
of the numerator and the denominator of the approximant. The fortunate choice
of the normalization condition depends on the particular situation.

A specific situation arises when the degree of the numerator (or of the
denominator) of the approximant is equal to zero. In this case the unfortunate
choice of the normalization condition results in the following:
the error approximant becomes zero or is not well-defined.
For $n=0$ it is expedient to choose condition~(42), as it was done in the
example given above. For $m=0$ (the case of the polynomial approximation)
it is usually expedient to choose condition~(43). Otherwise the
situation will be reduced to solving the system of equations~($27'$)
in the case described in~\S6 above.

Since the double precision regime of ES--1045 corresponds to 16 decimal
digits of mantissa in the computer representation of numbers,
while running computers of this type it makes sense to vary the
normalization condition only
in case the condition number exceeds $\delta\cdot 10^{16}$, where $\delta$
is the relative error of the obtained approximant. The value of the
condition number of the corresponding
system of linear algebraic equations is given by the PADE program
simultaneously with other computation results.

The theoretical error of the method is determined, to a considerable extent,
by the behavior of the approximant's denominator. It is
convenient for the analysis, by dividing the numerator and the denominator
of the fraction
by $b_0$ to equate $b_0$ to $1$. If the coefficients
$b_1,b_2,\dots ,b_m$ are small in comparison with $b_0=1$, which often
happens in computation practice, then the absolute error $\Delta (x)$
and its numerator $\Phi (x)=f(x)Q(x)-P(x)$ are of the same
order, so that the minimization of $\Phi (x)$ leads to
the minimization
of the error $\Delta (x)$, see~\S9 above. Note that the coefficients of
approximant~(45) to the function $\arctg x$ on the segment $[-1,1]$
are not small in comparison with $b_0$. For example, for $m=n=3$ the
coefficient $b_1$ is almost one and half times greater than the
coefficient $b_0$.
Thus, as shown in Table~1, the errors of the approximant to $\arctg x$
obtained by means of the PADE program are several times greater than the
errors of the best approximants.

Note that sometimes it is possible to improve the denominator of
the approximant or to reduce the condition number of the corresponding
system of equations by extending the segment $[A,B]$ where
the function is approximated. Such an effect is observed, for
example, when approximants to some hyperbolic
functions are calculated.

Note that the replacement of the standard subroutine DGELG for solving
system of linear algebraic equations by another subroutine of the same kind
(for example, by the DECOMP program from~\cite{41}) does not essentially
affect the quality of approximants obtained by means of
the PADE program.

One could seek the numerator and the denominator of the approximant in
the form
$$
\aligned
P&=\sum_{i=0}^na_iT_i,\\
Q&=\sum_{j=0}^mb_jT_j,
\endaligned
\tag50
$$
where $T_i$ are the Chebyshev polynomials. In this case the system of linear
equations determining the coefficients would be better
conditioned. But the calculation of the polynomials of the form~(50) by, for
example, the Chenshaw method, results in lengthening the computation time,
although it has a favorable effect upon the error of calculations,
see~\cite{47, Chapter IV,\S9}. The transformation of the polynomials
$P$ and $Q$ from the form~(50) into the standard form~(35) also
requires additional efforts.

In practice it is more convenient to use approximants
represented in the form~(1), (44), or~(45), and calculate the fraction's
numerator and denominator according the Horner scheme. In this case
the normalization $a_n=1$ or $b_m=1$ allows to reduce the number of
multiplications. Thus the PADE program gives coefficients of the
approximant in the two forms: with the condition $b_0=1$ and with
one of the conditions $a_n=1$ or $b_m=1$
no matter which one of the conditions~(41)--(43) is actually used
while solving the system of equations of type~(39) or~(40).

The PADE program (and the corresponding algorithm) can be easily modified,
for example, to take into account the case when some coefficients are fixed
beforehand.
One can vary the systems of equations under consideration by changing the
weight $w(x)$, the interval where the functions are approximated,
and the system of orthogonal polynomials. By a certain increase in complexity
of the system of equations~(40)
it is possible to minimize the norm of the numerator $\Phi (x)$
of the error function $\Delta (x)$ in the Hilbert space
$L_w^2$ (see~\S5 above).

The use of the PADE program does not require that
the approximated function be expanded into a series or a continued
fraction beforehand.
Equations~(39) or~(40) and the quadrature formula~(46) show that the PADE
program uses only the values of the approximated function
$f(x)$ at the interpolation points of the quadrature formula (which are
zeros of some Chebyshev polynomial).

On the segment $[-1,1]$ the linear Pad\'e--Chebyshev approximants give a
considerably smaller error than the classical Pad\'e approximants. For
example, the Pad\'e approximant of the form~(1) to the function $e^x$ for
$m=n=2$ has the absolute error $\Delta (1)=4\cdot 10^{-3}$ at the point
$x=1$, but the PADE program gives an approximant of the same form with
the absolute error $\Delta =1.9\cdot 10^{-4}$ (on the entire the
segment), i.e., the latter is 20 times smaller than the previous one.
The absolute error of the best approximant is $0.87\cdot 10^{-4}$.

\head\S11. The ``cross--multiplied''
linear Pad\'e--Chebyshev approximation scheme\endhead

As a rule, linear Pad\'e--Chebyshev approximants are constructed
according to the following scheme~\cite{45, 3, 11, 12}. Let the
approximated function be decomposed into the series in Chebyshev
polynomials
$$
f(x)=\sumpr_{i=0}^\infty c_iT_i(x)={1\over
2}c_0+c_1T_1(x)+c_2T_2(x)+\dots,
\tag51
$$
where the notation $\sumpr_{i=0}^m u_i$ means that the first term $u_0$
in the sum is replaced by $u_0/2$. The rational approximant is
looked for in the form
$$
R(x)={\sumpr\limits_{i=0}^na_iT_i(x)\over
\sumpr\limits_{j=0}^mb_jT_j(x)};\tag52
$$
the coefficients $b_j$ are determined by means of the system of linear
algebraic equations
$$
\sumpr_{j=0}^m b_j(c_{i+j}+c_{\vert i-j\vert})=0,\qquad
i=n+1,\dots,n+m,\tag53
$$
and the coefficients $a_i$ are determined by the equalities
$$
a_i={\frac 12}\sumpr_{j=0}^m b_j(c_{i+j}+c_{\vert i-j\vert})=0,\qquad
i=0,1,\dots,n.\tag54
$$
It is not difficult to verify that this algorithm must lead to the
same results as the algorithm described in~\S9
if the calculation errors are not taken into account.

The coefficients $c_k$ for $k=0,1,\dots,n+2m$, are present in~(53) and~(54),
i.e., it is necessary to have the first $n+2m+1$ terms of series~(51).
The coefficients $c_k$ are known, as a rule, only approximately. To determine
them one can take the truncated expansion of $f(x)$ into the series in
powers of $x$ (the Taylor series) and by means of the economization
procedure transform it into the form
$$
\sum_{i=0}^{n+2m} \tilde c_iT_i(x).\tag55
$$

\head\S12. Nonlinear Pad\'e--Chebyshev approximations\endhead

A rational function $R(x)$ of the form~(1) or~(52) is called a
{\it nonlinear Pad\'e--Chebyshev approximant} to the function $f(x)$
on the segment $[-1,1]$, if
$$
\int\limits_{-1}^1 \big(f(x)-R(x)\big)T_k(x)w(x)\,dx=0,\qquad
k=0,1,\dots,m+n,\tag56
$$
where $T_k(x)$ are the Chebyshev polynomials, $w(x)=1/\sqrt{1-x^2}$.
Che\-ney's theorem (see \S5 above) shows that the absolute error
function $\Delta (x)=f(x)-R(x)$ has alternation. Thus, there are
reasons to assume that the nonlinear Pad\'e--Chebyshev approximants are
close to the best ones in the sense of the absolute error.

In the paper~\cite{32} the following algorithm of computing the
coefficients of the approximant indicated above is given. Let the
approximated function $f(x)$ be expanded into series~(51) in 
Chebyshev polynomials. Determine the auxiliary quantities $\gamma_i$ from
the system of linear algebraic equations
$$
\sum_{j=0}^m\gamma_jc_{\vert k-j\vert}=0,\qquad
k=n+1,n+2,\dots,n+m,\tag57
$$
assuming that $\gamma_0=1$. The coefficients of the denominator in
expression~(52) are determined by the equalities
$$
b_j=\mu\sum_{i=0}^{m-j}\gamma_i\gamma_{i+j},
$$
where $\mu^{-1}=1/2\sum_{i=1}^n\gamma_i^2$; this implies
$b_0=2$. Finally, the coefficients of the numerator are determined by 
formula~(54). It is possible to solve system~(57) explicitly and to 
indicate the
formulas for computing the quantities $\gamma_i$. One can also estimate
explicitly the absolute error of the approximant. This
algorithm is described in detail in the book~\cite{33}; see also~\cite{11}.

In contrast to the linear Pad\'e--Chebyshev approximants, the nonlinear 
approximants of this type do not always exist, but it is possible to indicate 
explicitly verifiable conditions guaranteeing the 
existence of such approximants~\cite{33}. The nonlinear Pad\'e--Chebyshev
approximants (in comparison with the linear ones) have, as a rule, a somewhat
smaller absolute errors, but can have larger relative errors. 
Consider, as an example, the approximant of the form~(1) or~(52)
to the function $e^x$ on the segment $[-1,1]$ for $m=n=3$. In this case
the absolute error for a nonlinear Pad\'e--Chebyshev approximant 
is $\Delta =0.258\cdot 10^{-6}$, and the relative error,
$\delta =0.252\cdot 10^{-6}$; for the linear Pad\'e--Chebyshev approximant 
$\Delta =0.33\cdot 10^{-6}$ and $\delta =0.20\cdot 10^{-6}$.

\head\S13. Applications of the computer algebra system 
REDUCE to the construction of rational approximants\endhead

The computer algebra system REDUCE~\cite{48, 49} allows to handle
formulas at symbolic level and is a convenient tool for the
implementation of algorithms of computing rational approximants.
The use of this system allows to bypass 
the procedure of working out the algorithm of computing the approximated 
function if this function is presented in
analytical form or when either the Taylor series coefficients are
known or are determined analytically from a differential equation.
The round-off errors can be eli\-minated by using the exact  
arithmetic of rational numbers represented in the form of 
ratios of integers.

Within the framework of the REDUCE system, the program package for 
enhanced precision computations and construction of rational approximants 
is implemented; see, for example~\cite{8}. In particular, the algorithms
from~\S11 and~\S12 (which are similar to each other in structure) are
implemented, the approximated function being first expanded into the
power (Taylor) series,
$f=\sum_{k=0}^{\infty}f^{(k)}(0){x^k/k!}$, and then the truncated series 
$$
\sum_{k=0}^N f^{(k)}(0){x^k\over k!},\tag58
$$
consisting of the first $N+1$ terms of the Taylor series (the value $N$
is determined by the user) being transformed into a polynomial
of the form~(55) by means of the economization procedure.

The algorithms implemented by means of the REDUCE system allow to obtain 
approximants in the form~(1) or~(52), estimates of the absolute and the
relative error, and the error curves. The output includes the Fortran 
program of computing the corresponding approximant, the constants of 
rational arithmetic being transformed into the standard floating point
form. When computing the values of the obtained approximant, this 
approximant can be transformed into the form most convenient for the 
user. For example, one can calculate values of the numerator and the
denominator of the fraction of the form~(1) according to the Horner scheme,
and for the fraction of the form~(52), according to Clenshaw scheme, and
transform the rational expression into a continued fraction or a Jacobi
fraction as well.

The ALGOL-like input language of the REDUCE system and convenient tools for
solving problems of linear algebra guarantee simplicity and
compactness of the programs. For example, the length of the program for 
computing linear Pad\'e--Chebyshev approximants is sixty two lines.

\head\S14. The effect of error autocorrection for nonlinear
Pad\'e--Chebyshev approximations\endhead

Relations~(56) can be regarded as a system of equations for 
the coefficients of the approximant. Let the approximants
$R(x)=P(x)/Q(x)$ and $\widetilde{R}(x)=\widetilde{P}(x)/\widetilde{Q}(x)$,
where $P(x)$, $\widetilde{P}(x)$ are polynomials of degree $n$ and
$Q(x)$, $\widetilde{Q}(x)$ are polynomials of degree $m$, be obtained 
by approximate solving the indicated system of equations. Consider the
error approximant $\Delta P(x)/\Delta Q(x)$, where 
$\Delta P(x)=\widetilde{P}(x)-P(x)$, $\Delta Q(x)=\widetilde{Q}(x)-Q(x)$.
Substituting $R(x)$ and $\widetilde{R}(x)$ in~(56) and subtracting one of the 
obtained expressions from the other, we see that the following 
approximate equality holds:
$$
\int\limits_{-1}^1\bigg({\tilde P(x)\over\tilde Q(x)}-{P(x)\over 
Q(x)}\bigg)
T_k(x)w(x)\,dx\approx 0,\qquad k=0,1,\dots,m+n.
$$
This and equality~(29) imply the approximate equality
$$
\int\limits_{-1}^1\bigg({\Delta P(x)\over\Delta Q(x)}-{P(x)\over 
Q(x)}\bigg)
{\Delta Q\over\tilde Q}T_k(x)w(x)\,dx\approx \myone,\tag59
$$
where $k=0,1,\dots,m+n$, $w(x)=1/\sqrt{1-x^2}$. If the quantity 
$\Delta Q$ is relatively not small (this is connected with the fact that  
the system of equations~(57) is ill-conditioned), then, as follows from 
equality~(59), we can naturally expect that the error approximant 
is close to $P/Q$ and, consequently, to the approximated function~$f(x)$.

Due to the fact that the arithmetic system of rational numbers is used, the
software described in~\S13 allows to eliminate the round-off errors and to 
estimate the ``pure'' influence of errors in the approximated function 
on the coefficients of the nonlinear 
Pad\'e--Chebyshev approximant. In 
this case the effect of error autocorrection can be substantiated by a more
accurate reasoning which is valid both for nonlinear Pad\'e--Chebyshev
approximants and for linear ones, and even for the linear generalized 
Pad\'e approximants connected with different systems of orthogonal 
polynomials. This reasoning is analogous to Y.~L.~Luke's 
considerations~\cite{5} given  in~\S8 above.

Assume that the function $f(x)$ is expanded into series~(51) and that the
rational approximant $R(x)=P(x)/Q(x)$ is looked for in the form~(52).

Let $\Delta b_j$ be the errors in coefficients of the approximant's
denominator~$Q$. In the linear case these errors arise when solving the
system of equations~(53), and in the nonlinear case, when
solving the system of equations~(54). In both the cases the coefficients 
in the approximant's numerator are determined by equations~(54), 
whence we have
$$
\Delta a_i={\frac 12}\sumpr_{j=0}^m\Delta b_j(c_{i+j}+
c_{\vert i-j\vert}),\qquad i=0,1,\dots,n.\tag60
$$
This implies the following fact: the error approximant
$\Delta P/\Delta Q$ satisfies the relations
$$
\int\limits_{-1}^1\big(f(x)\Delta Q(x)-\Delta 
P(x)\big)T_i(x)w(x)\,dx=0,\qquad
i=0,1,\dots,n,\tag61
$$
which are analogous to relations~(39) defining the linear 
Pad\'e--Chebyshev approximants. Indeed, let us use the well-known 
multiplication formula for Chebyshev polynomials:
$$
T_i(x)T_j(x)=\frac 12\big[T_{i+j}(x)+T_{\vert i-j\vert}(x)\big],\tag62
$$
where $i$, $j$ are arbitrary indices; see, for example~\cite{11--13, 33}.
Taking~(62) into account, the quantity $f\Delta Q-\Delta P$ 
can be rewritten in the following way:
$$
\align
f\Delta Q-\Delta P&=\bigg(\sumpr_{j=0}^m\Delta b_jT_j\bigg)\bigg(
\sumpr_{i=0}^\infty c_iT_i\bigg)-\sumpr_{i=0}^n\Delta a_iT_i\\
&=\frac 12\sumpr_{i=0}^\infty\bigg[\sumpr_{j=0}^m\Delta b_j(c_{i+j}+
c_{\vert i-j\vert})\bigg]T_i-\sumpr_{i=0}^n\Delta a_iT_i.
\endalign
$$
This formula and~(60) imply that 
$$
f\Delta Q-\Delta P=O(T_{n+1}),
$$
i.e., in the expansion of the function $f\Delta Q-\Delta P$ into the
series in Chebyshev polynomials, the first $n+1$ terms are absent,
and the latter is equivalent to relations~(61) by virtue of the fact 
that the Chebyshev polynomials form an orthogonal system. When carrying out
actual computations, the coefficients $c_i$ are known
only approximately, and thus the equalities~(60), (61) are also 
satisfied approximately.

Consider the results of computer experiments
\footnote {At the author's request, computer calculations 
were carried out by A.~Ya.~Rodionov.}
that were performed by means of the software implemented
within the framework of the REDUCE system and briefly described 
in~\S13 above. We begin with the example considered in~\S10 above,
where the linear Pad\'e--Chebyshev approximant of the form~(44) to the
function $\cos{\pi\over 4}x$ was constructed on the segment
$[-1,1]$ for $m=2$, $n=3$. To construct the corresponding nonlinear 
Pad\'e--Chebyshev approximant, it is necessary to specify the value of the
parameter $N$ determining the number of 
terms in the truncated Taylor series~(58) of the approximated function. 
In this case the  calculation error is determined, in fact, by the 
parameter $N$.

The coefficients in approximants of the form~(44) which are obtained for
$N=15$ and $N=20$ (the nonlinear case) and the coefficients in the
error approximant are as follows
\footnote{Here we have in mind the coefficients of the expansions of the 
approxi\-mant's numerator and denominator in powers of $x$.}:

$$
\aligned
\tilde a_0&=0.4960471034987563,\\
a_0&=0.4960471027757504,\\
\tilde a_1&=-0.1451091945278387,\\
a_1&=-0.1451091928755736,\\
\tilde a_2&=10^{-2}\cdot 0.5482586543334515,\\
a_2&=10^{-2}\cdot 0.548258121085953,\\
\tilde a_3&=-10^{-4}\cdot 0.5205903601778259,\\
a_3&=-10^{-4}\cdot 0.5205902238186334,\\
\tilde b_0&=0.4960471034987759,\\
b_0&=0.4960471027757698,\\
\tilde b_1&=10^{-2}\cdot 0.7884201590727615,\\
b_1&=10^{-2}\cdot 0.7884203019999351,\\
\tilde b_2&=10^{-4}\cdot 0.4977097973870693,\\
b_2&=10^{-4}\cdot 0.4977100977750249,
\endaligned
\qquad
\aligned
\Delta a_0&=10^{-8}\cdot 0.07230059,\\
{}\\
\Delta a_1&=-10^{-8}\cdot 0.16522651,\\
{}\\
\Delta a_2&=-10^{-9}\cdot 0.42224856,\\
{}\\
\Delta a_3&=-10^{-10}\cdot 0.13635919,\\
{}\\
\Delta b_0&=10^{-8}\cdot 0.07230061,\\
{}\\
\Delta b_1&=-10^{-10}\cdot 0.1429272,\\
{}\\
\Delta b_2&=-10^{-10}\cdot 0.300388.
\endaligned
$$
Both the approximants have absolute errors $\Delta$ equal to 
$0.4\cdot 10^{-13}$
and the relative errors $\delta$ equal to $0.6\cdot 10^{-13}$, these 
values being close to the best possible.
The condition number of the system of equations~(57) in
both the cases is $0.4\cdot 10^8$. The denominator $\Delta Q$ of the error 
approximant is zero for $x=x_0\approx 0.70752\dots$; the
point $x_0$ is also close to the root of the numerator $\Delta P$ which for
$x=x_0$ is of order $10^{-8}$. Such a situation was 
considered in~\S10 above. Outside a small neighborhood of
the point $x_0$ the absolute and the relative errors have the
same order as in the ``linear case'' considered in~\S10.

Now consider the nonlinear Pad\'e--Chebyshev approximant of the form (44)
on the segment $[-1,1]$ to the function $\tg{\pi\over 4}x$ for 
$m=n=3$. In this case the Taylor series converges very slowly, and,
as the parameter $N$ increases, the values of coefficients of the rational
approximant 
undergo substantial (even in the first decimal digits) and intricate changes.
The situation is illustrated in Table~2,
where the following values are given: the absolute errors $\Delta$,
the absolute errors $\Delta_0$ of error approximants 
\footnote{A small neighborhood of the root of the polynomial
$\Delta Q$ is eliminated as before.}
(there the approximants are compared for $N=15$ and $N=20$, 
for $N=25$ and $N=35$,
for $N=40$ and $N=50$), and also the values of the condition number 
$\cond$ of the system of linear algebraic equations~(57). In this case
the relative errors coincide with the absolute ones. The best possible 
error is $\Delta_{\min}=0.83\cdot 10^{-17}$.

\centerline{ Table 2}

$$
\eightpoint
\matrix
\format\c&\quad\c&\quad\c&\quad\c&\quad\c&\quad\c&\quad\c\\
\text{ N }&15&20&25&35&40&50\cr
\cr
\cond&0.76\cdot 10^7&0.95\cdot 10^8&0.36\cdot 10^{10}&0.12\cdot 10^{12}&
0.11\cdot 10^{12}&0.11\cdot 10^{12}\cr
\Delta&0.13\cdot 10^{-4}&0.81\cdot 10^{-6}&0.13\cdot 10^{-7}&0.12\cdot 
10^{-10}&
0.75\cdot 10^{-12}&0.73\cdot 10^{-15}
\endmatrix
$$
$$
\eightpoint
\matrix
\format\c&\c&\c&\c\\
\Delta_0&\qquad\qquad 0.7\cdot 10^{-4}&\qquad
\quad\qquad\qquad 0.7\cdot 10^{-8}&\qquad\qquad\qquad\qquad 0.2\cdot 10^{-9}
\qquad
\endmatrix
$$

\newpage

\head\S15. Small deformations of approximated functions and
acceleration of convergence of series\endhead

Let a function $f(x)$ be expanded into the series in Chebyshev
polynomials, $f(x)=\sum_{i=0}^\infty c_iT_i$; consider a partial sum
$$
\hat f_N(x)=\sum_{i=o}^Nc_iT_i\tag63
$$
of this series. Using formula~(62), it is easy to verify that the 
linear Pad\'e--Chebyshev
approximant of the form~(1) or~(52) to the function $f(x)$ coincides
with 
the linear Pad\'e--Chebyshev approximant to polynomial~(63) for
$N=n+2m$, i.e., it depends only on the first $n+2m+1$ terms of the 
Fourier--Chebyshev series of the function $f(x)$; a similar result is
valid for the approximant of the form~(44) or~(45) to even or odd 
functions, respectively. Note that for $N=n+2m$ the polynomial
$\hat f_N$ 
is the result of application of the algorithm of linear (or nonlinear)
Pad\'e--Chebyshev approximation to $f(x)$, the exponents $m$ and $n$ 
being replaced by $0$ and $2m+n$.

The interesting effect mentioned in~\cite{9} consists in the fact that 
the error 
of the polynomial approximant $\hat f_{n+2m}$ depending on $n+2m+1$
parameters can exceed the error of the corresponding Pad\'e--Chebyshev
approximant of the form~(1) which depends only on $n+m+1$ parameters.
For example, consider an approximant of the form~(45) to the function
$\tg{\pi\over 4}x$ on the segment $[-1,1]$. For $m=n=3$ the linear 
Pad\'e--Chebyshev approximant to $\tg{\pi\over 4}x$ has the error of
order $10^{-17}$, and the corresponding polynomial approximant of the
form~(63) has the error of order $10^{-11}$. This polynomial of degree 19
\footnote{Odd functions are in question, and hence
$m=n=3$ in~(45) corresponds to $m=6$, $n=7$ in~(1).}
can be regarded as a result of deformation of the approximated function
$\tg{\pi\over 4}x$. This deformation does not affect the first 
twenty terms in the expansion of this function in Chebyshev polynomials 
and, consequently, does not affect the coefficients in the corresponding 
rational Pad\'e--Chebyshev approximant, but leads to a several orders 
increase of its error. Thus, a small deformation of the approximated 
function can
result in a sharp change in the order of error of a rational approximant.

Moreover the effect just mentioned means that the
algorithm extracts from polynomial~(63) 
additional information concerning the next components of the 
Fourier--Chebyshev
series. In other words, in this case the transition from Fourier--Chebyshev
series to Pad\'e--Chebyshev approximant accelerates
convergence of series. A similar effect of acceleration of
convergence of power series by passing to the classical Pad\'e approximant is
known (see~\cite{11, 14, 15}).

It is easy to see that the nonlinear Pad\'e--Chebyshev approximant
of the form~(1) to the function $f(x)$ depends only on the first
$m+n+1$ terms of the Fourier--Chebyshev series for $f(x)$, so that for such
approximants a more pronounced effect of the type indicated above 
takes place.

Since one can change the ``tail'' of the Fourier--Chebyshev series in
a quite arbitrary way without affecting the rational Pad\'e--Chebyshev
approximant, the effect of acceleration of convergence can take
place only for the series with an especially regular behavior (and
for the corresponding ``nice'' functions).

Note that the effect of error autocorrection indicates to the fact that 
the variation of an 
approximated function under deformations of a more general type may have
little effect on the rational approximant considered as a function 
(whereas the
coefficients of the approximant can have substantial changes). Accordingly,
while deforming the functions for which good rational approximation is 
possible, the approximant's error can rapidly increase.

There are interesting results distinguishing the classes of functions for 
which an efficient rational approximation is possible, for example, the 
classes
of functions which are approximated by rational fractions considerably
better (with a higher rate of convergence),
then by polynomials; see, in particular,~\cite{10, 50--52}. 
The reasoning given above
indicate that of a special interest are ``individual'' properties of 
functions which guarantee their effective
rational approximation.
There are reasons to suppose that solutions of certain functional and 
differential equations possess properties of this kind. Note that in 
papers~\cite{16, 37}, starting from the fact that elementary functions 
satisfy simple differential equations, it is shown that these 
functions are better approximated by rational fractions than by polynomial 
ones (we have in mind the best approximation); 
because of complicated calculations only the following cases were 
considered: the denominator of a rational approximant is a linear function
or (for even and odd functions) is a polynomial of degree~2.

\head\S16. Applications to computer calculation \endhead Ti computer calculation of function values is reduced in fact to
carrying out a finite set of arithmetic operations with
the argument and constants, i.e. to computing the value of a certain 
rational function.
Now we list some typical applications of methods for constructing
rational approximants. 
Often it happens that a function $f(x)$ is to be 
computed many times (for example, when solving numerically a differential 
equation) and with a given accuracy. In this case the
construction of a rational approximant to this function (with a given 
accuracy) often produces the most economic algorithm for computation of
values of $f(x)$. For example, if $f(x)$ is a complicated aggregate
of elementary and special functions every one which can be
calculated using the corresponding standard programs, then values of 
the function $f(x)$ can, of course, be computed 
by means of these programs. But such
an algorithm is often too slow and produces an unnecessary extra precision.

Standard computer programs for elementary and special functions, in their
turn, are based, as a rule, on rational approximants. Note that although the 
accuracy of rational and polynomial approximants to a given function is
the same, the computation of the
rational approximant usually requires a lesser number of 
operations, i.e., it is more speedy; see, for
example~\cite{1, 3, 12, 13, 24, 25, 31}.

The coefficients of rational approximants to basic elementary and special
functions can be found in reference handbooks; we note especially
the fundamental book~\cite{3}, see, also, for example,~\cite{12, 13}. But a
computer can have certain specific properties requiring 
algorithms and approximants (for effective standard programs of computing
functions) which are absent in
reference handbooks. In that case the construction programs for rational
approximants, including the PADE program described in~\S9 above (see 
also~\cite{1, 7, 9}), can be useful.

For example, decimal computers (including calculators) are
widely used at pre\-sent. The reason is that the use of decimal arithmetic 
system (instead of the standard binary one) enables the user to avoid a 
considerable 
loss of computing time needed for the transformation of numbers from the 
decimal representation to the binary one and vice versa. This is 
especially important if the amount of the 
input/output operations is relatively large; the latter situation is 
characteristic for calculations in the interactive mode. 
A method of computing elementary functions on decimal computers which 
uses the technique of rational approximants is described in the 
Appendix below. The main idea of this method consists in the fact
that the computation of values of various elementary functions, by means
of simple algorithms, is reduced to the computation of a rational
function of a fixed form. Roughly speaking, all basic elementary functions 
are calculated according to the same formula. Only the coefficients of the  
rational expression depend on a calculated function.

\head\S17. Nonlinear models and rational approximants\endhead

One of the main problems of mathematical modeling is to construct 
analytic formulas (models) that approximately describe the functional
dependence between different quantities according to given 
``experimental'' data
concerning the values of these quantities. In particular, let the set of
real numbers $x_1,\dots,x_{\nu}$ which are values of the ``independent''
variable $x$ be given, and for every value $x_i$ of this variable the value 
$y_i$ of the ``dependent'' variable $y$ be given. The problem is to
construct a function $y=F(x)$ such that the functional dependence
can be represented by an analytic formula of a certain form, and
the approximate equality $y_i\approx F(x_i)$ be valid for all
$i=1,2,\dots,\nu$, where the function $F(x)$ should take ``reasonable
values'' at points $x$ lying between the given points $x_i$.
In practice the values $y_i$ are usually given with errors.

As it was noted above, computer calculation of functions is finally 
reduced to computation of some rational functions. Thus in many cases it
is natural to construct an analytic model in the form of the rational
function~(1), where the degrees of the numerator and the denominator and
also the values of the coefficients are determined in the process of 
modeling, see~\cite{14}. Of course, in this case we have in mind only the
one-factor models. One can construct multi-factor models by using 
rational functions of several variables.

If we have a simple program of constructing rational approximants
to continuous functions defined on finite segments of the real line,
then we can reduce the construction of a model to constructing 
rational approximants to a continuous function (although in numerical 
analysis, as a rule, the goal is to reduce continuous problems to discrete
ones). The construction of a model is carried out in two steps. On the
first step a continuous function $f(x)$ such that $f(x_i)=y_i$ is 
constructed. A linear or a cubic spline (depending on the user's
choice) is used as $f(x)$. The function whose graph coincides with the 
polygonal
line consisting of segments of straight lines that connect the points
$(x_i,y_i)$ with the coordinates $x_i$, $y_i$ is the linear spline;
the cubic spline is described, for example, in~\cite{41}. 
On the second step the 
model is constructed by means of the PADE program. This approach guarantees
the regular behavior of the model on the entire range of the argument.

If there are reasons to assume that the initial data lie
on a sufficiently smooth and regular curve, then it is expedient to use
a cubic spline. And if there are reasons to assume that the initial data
contain considerable errors or deviations from theoretically admissible 
data, then it is expedient to use a linear spline: the behavior of
a cubic spline at intermediate points in this case will be irregular.

The method for constructing models described above was implemented 
(together with I.~A.~Andreeva) as the SPLINE--PADE program. This 
program prints out the graphs of splines and rational approximants 
(together with the initial data), and this facilitates the analysis 
of models.
Of course, while choosing and analyzing models, it is necessary to
take into account the theoretical requirements on the model which are 
connected with specific features of a particular problem.

\example{Example} Let the points $x_1,\dots,x_\nu$ be uniformly distributed on the
segment $[-{\pi\over 4},{\pi\over 4}]$, $x_1=-{\pi\over 4}$, $\nu=32$, 
$x_\nu={\pi\over 4}$, $y_i=\cos x_i$. The rational approximant of the 
form~(45)
to the linear spline for $m=n=2$ gives an approximant to $\cos x$ on
$[-{\pi\over 4},{\pi\over 4}]$ with the absolute error $\Delta=10^{-3}$. 
If a cubic spline is applied, then the absolute
error $\Delta$ is $0.35\cdot 10^{-6}$ in this case.\endexample

Other approaches to the construction of models in the form of rational 
functions can be found, for example, in~\cite{14}.

The above results connected with the effect of error autocorrection
show that similar models can have quite different 
coefficients.
Thus the coefficients of models of this kind are, generally 
speaking, unstable; and one should be very careful when trying to give
a substantial interpretation for these coefficients.

\head {APPENDIX\\}
A method of implementation of basic calculations on decimal computers\endhead

\subhead 1. Introduction\endsubhead
A large relative amount of input/output operations is a
characteristic feature of modern interactive computer systems. 
This results in a waste of computing time of systems with binary number
representation: numbers are transformed from the decimal representation
to the binary one and vice versa. Therefore, certain computers use decimal
arithmetic system. As a rule, the use of decimal arithmetic system leads to a
decrease in the rate of calculations and to additional memory 
requirements connected with specific coding of decimal numbers.
The decrease in the rate of calculations is due to the
fact that the implementation of decimal operations, as compared to that of
binary ones, is more complicated; moreover, the binary representation is 
more
convenient for implementing algorithms for calculating certain functions
then the decimal one. Since the performance rate of floating point 
arithmetic
operations and the rate of calculating elementary functions determine,
to a considerable extent, the rate of mathematical data performing,
the quality of the corresponding algorithms is, especially for cheap
personal systems, of great importance.

Here we consider methods of implementation of the floating point 
arithmetic system and of organizing computations for elementary functions. 
These methods are convenient to
use on decimal computers (this pertains both to the
software and hardware implementation). They guarantee a sufficient
economy of memory simultaneously with a relatively high performance rate
of calculations. Examples of effective software implementation of 
these methods are given in~\cite{1, 53}. The hardware implementation 
is described in the patent~\cite{55}. The methods under consideration 
are also of interest for octal and hexadecimal computers.

\subhead 2. Floating point arithmetic system\endsubhead
When carrying out arithmetic operations with floating point numbers,
the exponents of these numbers undergo only the operations of addition, 
subtraction, and comparison. Almost all computers have means for these
operations since they are necessary for the command and the 
address codes operations. This fact provides an
opportunity to use the binary representation for the 
exponents when implementing the floating point arithmetic system. Since 
exponents are integers lying in certain bounds, the transformation of
exponents from binary to decimal representations does not encounter
serious obstacles. The choice of an appropriate
algorithm depends on the structure of a computer and the method
of coding of decimal numbers. For the standard coding 8421, when each
decimal digit corresponds to a binary tetrad, it is possible to use the 
fact
that in this case the numbers from 0 to 9 have the same coding in the
binary and the decimal representations. 
Therefore the binary representation $x_2$ of a number
$x$ can be converted into the decimal representation $x_{10}$ by 
successively subtracting
(in the binary arithmetic system) the numbers from 0 to 9 from 
$x_2$ and forming the number $x_{10}$ from the sums of these numbers (in 
the
decimal arithmetic system). Similarly, a decimal integer can be converted
into a binary one.

Binary representation of exponents enables one to save
memory, and the combination of decimal operations 
with more rapid binary operations of addition type enhances
the performance rate. As a rule, the software implementation of
the floating point arithmetic system leads to the fact that
floating point operations take two orders as much time
when compared with fixed point operations.
The implementation described in~\cite{1} is
much more efficient: for seven decimal digit numbers, the transition from
the fixed point to the floating point regimes results in double
computing time for multiplication and division, and to reduction of the
rate of addition and subtraction by one decimal order.

\subhead 3. The design of computation for elementary functions\endsubhead
The calculation of values of each of the basic elementary functions (at the
reduction stage) is reduced to calculation of values of an odd function on
a symmetric (with respect to the origin) interval.
This odd function is approximated by a rational fraction of the form
$$
R(y)=y{a+by^2+cy^4\over \alpha+\beta y^2+y^4},\tag1
$$
where $y$ is the reduced argument, and the coefficients $a$, $b$,
$c$, $\alpha$, $\beta$ depend on the approximated function. Thus all
algorithms of computation for basic elementary functions have the common
block~(1), and this fact guarantees an economy of memory. This
block can be implemented both as a carefully devised part of software or
as a part of hardware; this
can enhance the performance rate. For the reduction algorithms
described below, the approximant of the form~(1) can guarantee 8--9 
accurate
decimal digits. Because of specific features of a particular computer and
the way the common block is implemented, it can be required that
expression~(1) be transformed into a certain form, for example, into the 
form
$$
R(y)=y{a+y^2(b+cy^2)\over \alpha+y^2(\beta +y^2)},
\tag"{($1'$)}"
$$
or into a Jacobi fraction of the form
$$
R(y)=y\Bigg(c+{\frac{\mu}{y^2+\nu+{\frac{\text{\it 
\ae}}{y^2+\lambda}}}}\Bigg).
\tag"{($1''$)}"
$$

The calculation of elementary functions with enhanced precision is
organized according to a similar scheme. The approximant of the form~(1) is
replaced by the expression
$$
R(y)=y{a+by^2+cy^4+dy^6\over \alpha+\beta y^2+\gamma y^4+y^6}\tag2
$$
which can be transformed into the form similar to ($1'$) or ($1''$),
i.e.,
$$
\align
R(y)&=y{a+y^2(b+y^2(c+dy^2))\over \alpha+y^2(\beta +y^2(\gamma+y^2))},
\tag"{($2'$)}"\\
R(y)&=y\Bigg(d+{\xi\over y^2+\eta+{\mu\over y^2+\nu +{\text{\it \ae}\over
y^2+\lambda}}}\Bigg).\tag"{($2''$)}"
\endalign
$$

The coefficients $a$, $b$, $c$, $d$, $\alpha$, $\beta$, $\gamma$,
$\xi$, $\eta$, $\mu$, $\nu$, {\it \ae}, $\lambda$ in formulas~($1'$),
($1''$), (2), ($2'$), ($2''$) are constants that depend on the 
approximated
function. The approximants of the form~(2), ($2'$) or ($2''$) 
guarantee 12--13 accurate decimal digits
\footnote{Of course, the values of the argument for which the loss
of precision is inevitable are an exception. For example, if
$x=1+\Delta x$, then $\ln x\approx\Delta x$, and the number of significant
digits of $\ln x$ is smaller than the number of significant
digits of the argument $x$ by the number of zeros after the decimal point 
in the number $\Delta x$.}.

The reduction algorithms are uniform; in particular, for calculations with
ordinary and enhanced precision the same reduction algorithms are used.
These algorithms are described in~section 4 below. The errors of
approximants and values of the coefficients in expressions~(1), (2) and
in their modifications are given below. These coefficients are either
taken from~\cite{3}, or calculated by means of the PADE program
described in~\S9 above.

\subhead 4. Algorithms\endsubhead
The relative, mean relative, absolute, and mean absolute errors
are denoted by $\delta$, $\bar\delta$, $\Delta$, $\bar\Delta$,
respectively.

\subsubhead {\rom {4.1.}\/} Calculation of logarithms\endsubsubhead
Let the argument $x>0$ have the form
$x=x_0\cdot 10^p$, where $0.1\leq x<1$, $p$ is an integer.
Suppose 
$$
y={x_0-{\sqrt{10}\over 10}\over x_0+{\sqrt{10}\over 10}},
$$ 
then we have
$$
x_0=10^{-{1\over 2}}\cdot{1+y\over 1-y},
$$ 
whence
$$
\lg x=p-{\frac 12}+\lg{1+y\over 1-y}.
$$ 
Substituting the approximant of the form
$R(y)$ with the best possible absolute error for the odd function
$$
\lg{1+y\over 1-y}={2\Arcth y\over\ln 10},
$$ 
we finally obtain $\lg (x)\approx p-1/2+R(y)$ for
$$
{1-\sqrt{10}\over 1+\sqrt{10}}\leq y<{\sqrt{10}-1\over\sqrt{10}+1}.
$$ 
For $0.1\leq x\leq 1$ and ordinary precision, $\Delta =0.23\cdot 10^{-8}$,
$\bar\Delta =0.14\cdot 10^{-8}$. For enhanced precision,
$\Delta =0.85\cdot 10^{-12}$, $\bar\Delta =0.53\cdot 10^{-12}$. 
It is impossible to minimize the relative error on the given interval
since this error is inevitable in a neighborhood of the point $x=1$.

The calculation of the natural logarithm is reduced to the case of the
decimal logarithm by means of the relation $\ln x=(\ln 10)(\lg x)$.

\subsubhead {\rom {4.2.}\/} Calculation of exponentials\endsubsubhead
Consider a nonstandard (at the
first sight) algorithm of reduction of the function $10^x$, which,
nevertheless,
is dual to the algorithm of reduction of $\lg x$ described above.
Represent the argument $x$ in the form $x=y+p$, where $-1<y<1$, $p$ is
an integer (for example, $-12.85=-12+(-0.85))$. Then
$$
10^x\approx 10^p{1+R(y)\over 1-R(y)},
$$ 
where $R(y)$ is the approximant of the
form~(1) or~(2) on the interval $[-1,1]$ to the odd function
$\th (x\ln 10/2)$. For $-1\leq x\leq 1$ and ordinary precision,
$\delta =0.23\cdot 10^{-8}$, $\bar\delta =0.6\cdot 10^{-9}$,
$\Delta =0.23\cdot 10^{-7}$, $\bar\Delta =0.17\cdot 10^{-8}$. For
enhanced precision,
$\delta =0.17\cdot 10^{-13}$, $\bar\delta =0.3\cdot 10^{-14}$,
$\Delta =0.16\cdot 10^{-12}$, $\bar\Delta =0.85\cdot 10^{-14}$.

For calculation of the functions $e^x$ and $x^y$, the relations
$e^x=10^{x\lg e}$ and $x^y=10^{y\lg x}$ are used.

\subsubhead {\rom {4.3.}\/} Calculation of $\sin x$ and $\cos x$\endsubsubhead
Denote by $R(x)$ the approximant of the
form~(1) or~(2) with the best possible relative error to the function 
$\sin x$ on
the segment $[-\pi/2,\pi/2]$. Since the function $\sin x$ is odd,
it is sufficient to consider the case $x>0$. Denote by $\{a\}$ the
fractional part of a positive number $a$; for example, $\{12.08\}=0.08$.
Set $y=\{x/2\pi\}\cdot 2\pi$, then $0\leq y<2\pi$. If
$y\leq\pi/2$, then we set $\sin x\approx R(y)$; if
$\pi/2<y\leq 3\pi/2$, then $\sin x\approx R(z)$, where $z=\pi-y$,
and if $3\pi/2<y<2\pi$, then $\sin x\approx R(z)$, where $z=y-2\pi$.
For $-\pi/2\leq x\leq\pi/2$ and ordinary precision,
$\delta =\Delta =0.53\cdot 10^{-8}$, $\bar\delta =0.34\cdot 10^{-8}$,
$\bar\Delta =0.21\cdot 10^{-8}$. For enhanced precision,
$\delta =\Delta =0.56\cdot 10^{-13}$, $\bar\delta =0.32\cdot 10^{-13}$,
$\bar\Delta =0.2\cdot 10^{-13}$. The calculation of $\cos x$ is
reduced to the calculation of $\sin x$ by means of the relation
$\cos x=\sin (\pi/2-x)$.

\subsubhead {\rom {4.4.}\/} Calculation of $\tg x$\endsubsubhead
Let $R(x)$ be the approximant of the
form~(1) or~(2) with the best possible relative error to the function 
$\tg x$
on the segment $[-\pi/4,\pi/4]$, the approximant~(2) 
satisfying the additional condition $d=0$. The algorithm of reduction
is quite similar to the algorithm for $\sin x$ given above. For
$x>0$ set $y=\{x/\pi\}\pi$; in this case $0\leq y<\pi$.
Hence $\tg x\approx R(y)$ for $y\leq\pi/4$; for $\pi/4<y\leq 3\pi/4$
we have
$\tg x\approx 1/R(z)$, where $z=\pi/2-y$; finally, for
$3\pi/4<y<\pi$ we have $\tg x\approx R(z)$, where $z=y-\pi$. 
For $x<0$ we use the
relation $\tg(-x)=-\tg x$. For $-\pi/4\leq x\leq\pi/4$ and ordinary
precision
$\Delta =\delta =0.22\cdot 10^{-10}$, $\bar\Delta =0.63\cdot 10^{-11}$,
$\bar\delta =0.13\cdot 10^{-10}$. For enhanced precision,
$\Delta =\delta =0.26\cdot 10^{-13}$, $\bar\delta =0.15\cdot 10^{-13}$,
$\bar\Delta =0.67\cdot 10^{-14}$. The 
algorithm for calculating $\tg x$ described above has essential 
advantages in accuracy and speed as compared with the algorithm using
the relation $\tg x=\sin x/\cos x$ and the algorithms for calculating
$\sin x$ and $\cos x$.

\subsubhead {\rom {4.5.}\/} Calculation of $\arctg x$\endsubsubhead
Let $R(x)$ be the approximant of the form (1)
or (2) with the best possible relative error to the function $\arctg x$ for
$\vert x\vert\leq\tg\pi/8=\sqrt{2}-1$. The reduction is standard:
if $0\leq x<\sqrt{2}-1$, then $\arctg x\approx R(x)$; if $\sqrt{2}-1\leq 
x<1$,
then $\arctg x\approx\pi/4-R(y)$, where $y=(1-x)/(1+x)$; if
$x>1$, then $\arctg x\approx\pi/2-R(1/x)$; for $x<0$ the relation
$\arctg (-x)=-\arctg x$ is used. For $\vert x\vert\leq\sqrt{2}-1$ and 
ordinary
precision,
$\delta =0.29\cdot 10^{-9}$, $\bar\delta =0.18\cdot 10^{-9}$,
$\Delta =0.11\cdot 10^{-9}$, $\bar\Delta =0.37\cdot 10^{-10}$.
For enhanced precision,
$\Delta =0.11\cdot 10^{-13}$, $\delta =0.29\cdot 10^{-13}$,
$\bar\delta =0.18\cdot 10^{-13}$, $\bar\Delta =0.36\cdot 10^{-14}$.

\subsubhead {\rom {4.6.}\/} Calculation of $\arcsin x$\endsubsubhead
Let $R(x)$ be the approximant of the
form (1) or (2) with the best possible relative error to the function
$\arcsin{x}$ on the interval $[-1/2,1/2]$. Since the function
is odd, it is sufficient to consider the case $x>0$. If
$0<x\leq 1/2$, then $\arcsin x\approx R(x)$; if $1/2<x\leq 1$,
then $\arcsin x\approx\pi/2-2R(y)$, where $y=\sqrt{(1-x)/2}$.
For $-1/2\leq x\leq 1/2$ and ordinary precision,
$\delta =0.25\cdot 10^{-8}$, $\bar\delta =0.16\cdot 10^{-8}$,
$\Delta =0.13\cdot 10^{-8}$, $\bar\Delta =0.41\cdot 10^{-9}$; for
enhanced precision,
$\delta =0.82\cdot 10^{-12}$, $\bar\delta =0.52\cdot 10^{-12}$,
$\Delta =0.43\cdot 10^{-12}$, $\bar\Delta =0.13\cdot 10^{-12}$.

{\rom {4.7.}\/} {\it The reduction algorithms} 
for  $\lg x$, $\arcsin x$, and $\arctg x$
desc\-ri\-b\-ed above are taken from~\cite{3}. The reduction algorithms for 
$\sin x$ and $\tg x$ were proposed by the author and R.~M.~Borisyuk~\cite{53}.
Of course, in particular cases the general scheme is supplemented by
special ruses. For example, $\sin x$, $\arcsin x$, and $\arctg x$ are
replaced by $x$ for small values of the argument, and so on.

\subhead 5. Coefficients\endsubhead
For every function the coefficients of approximants that
 are used while computing values of this function are indicated below
(see Table 3 and Table~4).
For every function the coefficients of approximants~(1), ($1'$) and
($1"$) are listed according to the following order: $a$, $b$, $c$,
$\alpha$, $\beta$, $\gamma$, $d$, $\xi$, $\eta$, $\mu$, $\lambda$,
$\nu$, {\it\ae\/}; the coefficients of approximants~(2), ($2'$) and
($2"$) are listed according to the following order: $a$, $b$, $c$,
$d$, $\alpha$, $\beta$, $\gamma$, $\xi$, $\eta$, $\mu$, $\lambda$,
$\nu$, {\it\ae\/}; mantissas (significands) are separated from  exponents
by the letter $D$. The accuracy of
the coefficients (16 decimal digits of
the mantissa) is, of course, excessive.

\subhead 6. Analysis of the algorithms\endsubhead
It is easy to see that the algorithms of
calculating trigonometric and inverse trigonometrical functions do not
depend on on the arithmetic system of the computer. On the contrary, 
while implementing the computing algorithms for exponentials,
logarithms and functions that are expressed through them (hyperbolic and
inverse hyperbolic functions
\footnote{Note that it is also convenient to use the common
block of type (1) or (2) while calculating hyperbolic and inverse
hyperbolic functions.}, 
$x^y$) the binary arithmetic system has an essential
advantage over the decimal one. For example, for binary arithmetic system
the computation of the logarithm is reduced to finding an approximant on
the segment $[1/2,1]$ (and not on the segment $[1/10,1]$); since
1/2 is much closer to zero than 1/10, this implies that the approximation
rate increases considerably. While computing $\ln x$
according to the scheme described above on a binary computer,
the approximant of the form~(1)
which depending on five parameters can be replaced by a more exact 
approximant
(on a smaller segment) which depending only on three parameters. A similar
situation arises while calculating an exponential. 
But the 
use of the decimal arithmetic system leads to a certain equilibrium
between the difficulty of computing logarithmic and exponential functions,
on one hand, and trigonometric functions, on the other.
Thus in this case the use of a separate common block of the form~(1) or~(2)
is justified.

\subhead 7. Implementation of algorithms for 
calculating elementary functions\endsubhead
\newline For the software implementation it is expedient to use 
representations~($1''$) and~($2''$) for rational
approximants in the form of Jacobi fractions; this allows to minimize 
the number of arithmetic operations. The rate of computation of 
functions can be increased by implementing the calculation of Jacobi 
fractions mentioned above by means of the fixed point arithmetic system as 
described in~\cite{1}.

A method of hardware implementation for algorithms under consideration
is described in the patent~\cite{55}.
In this case it is expedient
to use representations~($1'$) and~($2'$) for rational approximants and
to carry out computations of the fraction numerator and denominator 
in parallel. 
For example, when computing expression~($2'$), the value
$y^2$ being computed beforehand, it is possible to use the 
summator to compute 
$\gamma +y^2$ and the multiplier to compute
$d\cdot y^2$ simultaneously. Then $y^2$ is multiplied by $(\gamma+y^2)$
and simultaneously the quantity $c$ is added to $d\cdot y^2$, and so on.
Under such an implementation, additional hardware requirements
are minimal since almost all computers have a summator and a
multiplier.

\newpage

\centerline{Table 3. Ordinary precision}
$$
{\eightpoint
\matrix
\format\c&\c&\c&\quad\c&\c&\c\\
{}&\qquad\qquad\lg x&{}&{}&\qquad\qquad 10^x&{}\cr
{}&0.3187822082024000\DD&\quad 01&\qquad{}&0.4184196402707361\DD&\quad 02\cr
-&0.2655808794660000\DD&\quad 01&\qquad{}&0.6132751585841820\DD&\quad 01\\
{}&0.2668632700470000\DD&\quad 00&\qquad{}&0.7526525036394230\DD&-01\\
{}&0.3670115625115000\DD&\quad 01&\qquad{}&0.3634346820241857\DD&\quad 02\\
-&0.4280973292830000\DD&\quad 01&\qquad{}&0.2138428615920360\DD&\quad 02\\
{}&\qquad\qquad ***&{}&{}&\qquad\qquad ***&{}\\
{}&0.2668632700470000\DD&\quad 00&\qquad &0.7526525036394230\DD&-01\\
-&0.1513374262751513\DD&\quad 01&\qquad &0.4523257934215174\DD&\quad 01\\
-&0.1459257686092834\DD&\quad 01&\qquad &0.8645663007168292\DD&\quad 01\\
-&0.2821715606737166\DD&\quad 01&\qquad &0.1273862315203531\DD&\quad 02\\
-&0.4474945619843138\DD&\quad 00&\qquad -&0.7379037474539062\DD&\quad 02\\
{}&\qquad\qquad\sin x&{}&{}&\qquad\qquad\tg x&{}\\
{}&0.2051458702138878\DD&\quad 04&\qquad &0.6260411195474330\DD&\quad 02\\
-&0.2731535822018325\DD&\quad 03&\qquad -&0.6971684006294421\DD&\quad 01\\
{}&0.6635500992122553\DD&\quad 01&\qquad &0.6730910258759150\DD&-01\\
{}&0.2051458712958973\DD&\quad 04&\qquad &0.6260411195336056\DD&\quad 02\\
{}&0.6875599504020228\DD&\quad 02&\qquad -&0.2783972122004270\DD&\quad 02\\
%%{}&\qquad\qquad ***&{}&{}&\qquad\qquad ***&{}\\
%%{}&0.6635500992122553\DD&\quad 01&\qquad &0.6730910258759150\DD&-01\\
%%-&0.7293840555054680\DD&\quad 03&\qquad ÿ17#58461	8D&\quad 01\
{}&\qquad\qquad ***&{}&{}&\qquad\qquad ***&{}\\
{}&0.6635500992122553\DD&\quad 01&\qquad &0.6730910258759150\DD&-01\\
-&0.7293840555054680\DD&\quad 03&\qquad -&0.5097817354684619\DD&\quad 01\\
{}&0.1585035693573942\DD&\quad 02&\qquad -&0.1145397751592885\DD&\quad 02\\
{}&0.5290563810446287\DD&\quad 02&\qquad -&0.1638574370411386\DD&\quad 02\\
{}&0.1212885465090180\DD&\quad 04&\qquad -&0.1250778280153322\DD&\quad 03\\
{}&\qquad\qquad\arctg x&{}&{}&\qquad\qquad\arcsin x&{}\\
{}&0.4482500977985320\DD&\quad 01&\qquad &0.5603629044813127\DD&\quad 01\\
{}&0.3372473379182700\DD&\quad 01&\qquad -&0.4614530946664500\DD&\quad 01\\
{}&0.2742666270116000\DD&\quad 00&\qquad &0.4955994747873100\DD&\quad 00\\
{}&0.4482500979270910\DD&\quad 01&\qquad &0.5603629030606043\DD&\quad 01\\
{}&0.4866639968788300\DD&\quad 01&\qquad -&0.5548466599346680\DD&\quad 01\\
{}&\qquad\qquad ***&{}&{}&\qquad\qquad ***&{}\\
{}&0.2742666270116000\DD&\quad 00&\qquad &0.4955994747873100\DD&\quad 00\\
{}&0.2037716450063295\DD&\quad 01&\qquad -&0.1864713814153353\DD&\quad 01\\
{}&0.1596444173439070\DD&\quad 01&\qquad -&0.1515767952641660\DD&\quad 01\\
{}&0.3270195795349230\DD&\quad 01&\qquad -&0.4032698646705020\DD&\quad 01\\
-&0.7381840442193144\DD&\quad 00&\qquad -&0.5090063407308178\DD&\quad 00
\endmatrix}
$$

\newpage

\centerline{Table 4. Enhanced precision}

$$
{\eightpoint
\matrix 
\format\c&\c&\c&\quad\c&\quad\c&\c\\
{}&\qquad\qquad\lg x&{}&{}&\qquad\qquad 10^x&{}\cr
-&$0.8625170319686105$\DD&\quad 01&\qquad{}&0.2416631060448244\DD&\quad 04\cr
{}&0.1170031513942458\DD&\quad 02&\qquad{}&0.4101315802439533\DD&\quad 03\cr
-&0.3932918863942010\DD&\quad 01&\qquad{}&0.1179791485292238\DD&\quad 02\cr
{}&0.1960631750250080\DD&\quad 00&\qquad{}&0.4067164397089984\DD&-01\cr
-&0.9930094301066197\DD&\quad 01&\qquad{}&0.2099059068697363\DD&\quad 04\cr
{}&0.1678051701465279\DD&\quad 02&\qquad{}&0.1283652206816926\DD&\quad 04\cr
-&0.8135425915212830\DD&\quad 01&\qquad{}&0.8569022438348670\DD&\quad 02\cr
{}&\qquad\qquad ***&{}&{}&\qquad\qquad ***&{}\cr
{}&0.1960631750250080\DD&\quad 00&\qquad &0.4067164397089984\DD&-01\cr
-&0.2337861428824651\DD&\quad 01&\qquad &0.8312752555010688\DD&\quad 01\cr
-&0.4538004104289327\DD&\quad 01&\qquad &0.4263308628906700\DD&\quad 02\cr
-&0.2401159236894577\DD&\quad 01&\qquad -&0.8324501520959055\DD&\quad 03\cr
-&0.1263136819683083\DD&\quad 01&\qquad &0.1184109379926876\DD&\quad 02\cr
-&0.2334284991240421\DD&\quad 01&\qquad &0.3121604429515093\DD&\quad 02\cr
-&0.9196001135281050\DD&\quad -01&\qquad -&0.8918843336784789\DD&\quad 02\cr
\allowmathbreak
{}&\qquad\qquad\sin x&{}&{}&\qquad\qquad\tg x&{}\cr
{}&0.5896178692831105\DD&\quad 06&\qquad -&0.1025905306253198\DD&\quad 05\cr
-&0.8390788060005352\DD&\quad 05&\qquad &0.1244879443625433\DD&\quad 04\cr
{}&0.2686819889613831\DD&\quad 04&\qquad -&0.2085510846487616\DD&\quad 02\cr
-&0.2413382332334080\DD&\quad 02&\qquad &0.0&{}\cr
{}&0.5896178692830811\DD&\quad 06&\qquad -&0.1025905306253223\DD&\quad 05\cr
{}&0.1436176428157609\DD&\quad 05&\qquad &0.4664563797829625\DD&\quad 04\cr
{}&0.1669650188299142\DD&\quad 03&\qquad -&0.2078359665454338\DD&\quad 03\cr
{}&\qquad\qquad ***&{}&{}&\qquad\qquad ***&{}\cr
-&0.2413382332334080\DD&\quad 02&\qquad &0.0&{}\cr
{}&0.6716324155233251\DD&\quad 04&\qquad -&0.2085510846487616\DD&\quad 02\cr
{}&0.1278518975583208\DD&\quad 03&\qquad -&0.1481441435217348\DD&\quad 03\cr
{}&0.7154609949974411\DD&\quad 04&\qquad -&0.4670350587972472\DD&\quad 04\cr
{}&0.4298163104533473\DD&\quad 02&\qquad -&0.1340714587130577\DD&\quad 02\cr
-&0.3868509773741323\DD&\quad 01&\qquad -&0.4628467715239315\DD&\quad 02\cr
{}&0.2372742417389960\DD&\quad 04&\qquad -&0.1286250294831805\DD&\quad 03\cr
{}&\qquad\qquad\arctg x&{}&{}&\qquad\qquad\arcsin x&{}\cr
{}&0.1360093213361806\DD&\quad 02&\qquad -&0.1650992722517542\DD&\quad 02\cr
{}&0.1696003978718046\DD&\quad 02&\qquad &0.2211306741883449\DD&\quad 02\cr
{}&0.5018932379116295\DD&\quad 01&\qquad -&0.7427989795595475\DD&\quad 01\cr
{}&0.2013116125542811\DD&\quad 00&\qquad &0.4054985920387412\DD&\quad 00\cr
{}&0.1360093213361844\DD&\quad 02&\qquad -&0.1650992722516183\DD&\quad 02\cr
{}&0.2149368383148177\DD&\quad 02&\qquad &0.2486472195166319\DD&\quad 02\cr
{}&0.9463307253423236\DD&\quad 01&\qquad -&0.1033386531178065\DD&\quad 02\cr
{}&\qquad\qquad ***&{}&{}&\qquad\qquad ***&{}\cr
{}&0.2013116125542811\DD&\quad 00&\qquad &0.4054985920387412\DD&\quad 00\cr
{}&0.3113858735833039\DD&\quad 01&\qquad -&0.3237621961350434\DD&\quad 01\cr
{}&0.5406247281512437\DD&\quad 01&\qquad -&0.6618033809642488\DD&\quad 01\cr
-&0.3928353166591151\DD&\quad 01&\qquad -&0.2758376635378124\DD&\quad 01\cr
{}&0.1338761180200926\DD&\quad 01&\qquad -&0.1288186871451920\DD&\quad 01\cr
{}&0.2718298791709874\DD&\quad 01&\qquad -&0.2427644630686242\DD&\quad 01\cr
-&0.1505853445310237\DD&\quad 00&\qquad -&0.9565986684443910\DD&-01
\endmatrix}
$$

\eject

\bigskip

E-mail: litvinov\@islc.msk.su

glitvinov\@mail.ru

\enddocument